\documentclass[aos,seceqn,nameyear,rotating,dvips]{arximspdf}
\usepackage{multirow}
\usepackage{dcolumn}
\usepackage{graphicx}

\doi{10.1214/09-AOS784}
\volume{38}
\issue{4}
\pubyear{2010}
\firstpage{2351}
\lastpage{2387}

\makeatletter

\newcolumntype{d}[1]{D{.}{.}{#1}}

\newtheorem{theorem}{Theorem}[section]
\newproclaim{definition}{Definition}

\newproclaim{Remark}{Remark}

\newtheorem{Proposition}{Proposition}
\newtheorem{Lemma}{Lemma}

\makeatother

\begin{document}
\begin{frontmatter}

\title{Sieve estimation of constant and time-varying
coefficients in nonlinear ordinary differential equation models by
considering both numerical error and measurement error}
\runtitle{Sieve estimation in ordinary differential equation models}

\begin{aug}
\author[A]{\fnms{Hongqi} \snm{Xue}\ead[label=e1]{Hongqi\_Xue@urmc.rochester.edu}},
\author[A]{\fnms{Hongyu} \snm{Miao}\ead[label=e2]{Hongyu\_Miao@urmc.rochester.edu}} and
\author[A]{\fnms{Hulin} \snm{Wu}\corref{}\thanksref{T1}\ead[label=e3]{hwu@bst.rochester.edu}}
\runauthor{H. Xue, H. Miao and H. Wu}
\affiliation{University of Rochester}
\address[A]{
Department of Biostatistics\\
\quad and Computational Biology\\
University of Rochester School\\
\quad of Medicine and Dentistry\\
601 Elmwood Avenue, Box 630\\
Rochester, New York 14642 \\
USA\\
\printead{e1}\\
\phantom{E-mail: }\printead*{e2}\\
\phantom{E-mail: }\printead*{e3}} 
\end{aug}

\thankstext{T1}{Supported in part by NIAID/NIH Grants AI055290,
AI50020, AI078498, AI078842 and
AI087135.}

\received{\smonth{8} \syear{2009}}
\revised{\smonth{12} \syear{2009}}

%
\begin{abstract}
This article considers estimation of constant and time-varying
coefficients in nonlinear ordinary differential equation (ODE)
models where analytic closed-form solutions are not available. The
numerical solution-based nonlinear least squares (NLS) estimator is
investigated in this study. A numerical algorithm such as the Runge--Kutta
method is used to approximate the ODE solution. The asymptotic
properties are established for the proposed estimators
considering both numerical error and measurement error. The
B-spline is used to approximate the time-varying
coefficients, and the corresponding asymptotic theories in this case
are investigated under the framework of the sieve approach. Our
results show that if the maximum\vspace*{1pt} step size of the $p$-order
numerical algorithm goes to zero at a rate faster than
$n^{-1/(p\wedge4)}$, the numerical error is negligible compared to
the measurement error. This result provides a theoretical guidance in
selection of the step size for numerical evaluations of ODEs.
Moreover, we have shown that the numerical solution-based NLS
estimator and the sieve NLS estimator are strongly consistent. The
sieve estimator of constant parameters is asymptotically normal with
the same asymptotic co-variance as that of the case where the true
ODE solution is exactly known, while the estimator of the
time-varying parameter has the optimal convergence rate under some
regularity conditions. The theoretical results are also developed
for the case when the step size of the ODE numerical solver does not
go to zero fast enough or the numerical error is comparable to the
measurement error. We illustrate our approach with both simulation
studies and clinical data on HIV viral dynamics.
\end{abstract}

%

\end{frontmatter}

\section{Introduction}
Ordinary differential equations (ODE) are widely used to model dynamic
processes in many scientific
fields such as engineering, physics, econometrics and biomedical
sciences. In particular, new
biotechnologies allow scientists to use ODE models to more accurately
describe biological processes
such as genetic regulatory networks, tumor cell kinetics, epidemics and
viral dynamics of infectious
diseases [\citet{r11}, \citet{r30}, \citet{r40},
\citet{r12}, Anderson and May (\citeyear{r2}), \citet{r7}, \citet{r42}]. The
mathematical modeling approach has made a great impact on these
scientific fields over the past
decades. For instance, ODE models have been used to quantify HIV viral
dynamics which resulted in
important scientific findings [\citet{r22}, \citet{r70},
Perelson et al. (\citeyear{r46}, \citeyear{r44})]. Comprehensive reviews
of the application of ODE models in HIV dynamics can be found in
\citet{r45}, \citet{r42}, \citet{r64} and \citet{r73}.

Although differential equation models have been widely used in
scientific research, very little
statistical research has been dedicated to parameter estimation and
inference for differential
equation models. The existing statistical methods for ODE models
include the nonlinear least squares
method [\citet{r4}, \citet{r67}, \citet{r5},
\citet{r34}], the smoothing-based techniques [\citet{r63},
\citet{r68}, \citet{r10}, \citet{r35}, \citet{r8}],
the principal differential analysis (PDA) [\citet{r52}, \citet{r21},
\citet{r50}, \citet{r53}, Varziri et
al. (\citeyear{r69})] and the Bayesian approaches [\citet{r51}, \citet{r28},
\citet{r15}]. However, very few of these
publications rigorously address the theoretical issues and study the
asymptotic properties of
the proposed estimators when both measurement error and numerical error
are significant. In this paper, we intend to investigate statistical
estimation methods for both constant and time-varying parameters in ODE
models and study the asymptotic properties of the proposed estimators
under the framework of the sieve approach.

Denote a general set of ODE models containing only constant
parameters as
%
%
\begin{equation}\label{Const_ODE}
\cases{ \dfrac{d {\mathbf X}(t)}{d t}={\mathbf F}\{t,
\mathbf{X}(t),{\bolds\beta}\},&\quad $\forall t\in[t_0,T]$,\vspace*{2pt}\cr
{\mathbf X}(t_0)={\mathbf X}_0, & }
\end{equation}
and denote a general set of ODE models with both constant and
time-varying parameters as
%
%
\begin{equation}\label{Time_ODE}
\cases{\dfrac{d {\mathbf X}(t)}{dt}={\mathbf F}\{t, {\mathbf
X}(t),\bolds{\beta},\eta(t)\},&\quad $\forall t\in[t_0,T]$,\vspace*{2pt}\cr
{\mathbf X}(t_0)={\mathbf X}_0, & }
\end{equation}
where ${\mathbf X}(t)=\{X_1(t),\ldots,X_K(t)\}^T$ is a $K$-dimensional
state variable vector,
$\bolds{\beta}$ is a $d$-dimensional vector of unknown constant
parameters with true value
$\bolds{\beta}_0$, $\eta(t)$ is an unknown time-varying parameter with
true value $\eta_0(t)$ (here we
only consider a single time-varying parameter, the proposed methodology
can be extended to include
multiple time-varying parameters although it is tedious and cumbersome
in notation), ${\mathbf F}(\cdot)=\{F_1(\cdot),\ldots,F_K(\cdot)\}^T$ is a vector
of differentiable
functions whose forms are known and
${\mathbf X}(t_0)={\mathbf X}_0$ is the initial value. Equations (\ref
{Const_ODE}) and (\ref{Time_ODE}) are
called state equations. Obviously, equation (\ref{Const_ODE}) is a
special case of (\ref{Time_ODE}).
The function ${\mathbf F}(t,\mathbf{X},{\bolds\beta})$ in
(\ref{Const_ODE}) or ${\mathbf F}(t,{\mathbf X},\bolds{\beta},\eta)$ in
(\ref{Time_ODE}) is assumed to fulfil
the Lipschitz assumption to
$\mathbf{X}$ [with the Lipschitz constant independent of the unknown
parameters $\bolds{\beta}$ and
$\eta(\cdot)$] ensuring existence and uniqueness of the solutions to
(\ref{Const_ODE}) and
(\ref{Time_ODE}) [see \citet{r18} and
\citet{r38}]. Let
${\mathbf X}(t,\bolds{\beta})$ and ${\mathbf X}(t,\bolds{\beta
},\eta(t))$ denote
the true solutions to
(\ref{Const_ODE}) and (\ref{Time_ODE}) for given $\bolds
{\beta}$
and $\eta(\cdot)$, respectively. We
usually use notation ${\mathbf X}(t)$ to denote ${\mathbf X}(t,\bolds
{\beta
}_0)$ or ${\mathbf X}(t,\bolds{\beta}_0,\eta_0(t))$ in this article.
Our objective is to
estimate the unknown parameters
$\bolds{\beta}$ and $\eta(\cdot)$ based on the measurements of the state
variables, ${\mathbf X}(t)$ or
their functions.

If a closed-form solution to (\ref{Const_ODE}) or (\ref
{Time_ODE}) is available, the standard
statistical approaches for nonlinear regression or time-varying
coefficient regression models can be
used to estimate unknown parameters. In practice, (\ref
{Const_ODE}) and (\ref{Time_ODE})
usually do not have closed-form solutions for a nonlinear ${\mathbf F}$.
In this case, numerical methods such as the Runge--Kutta algorithm
[\citet{r55}, \citet{r32a}] have to be used to approximate the solution
of the ODEs for a given set of parameter values and initial conditions.
Consequently, the nonlinear least squares (NLS) principle (minimizing
the residual sum of squares of the differences between the experimental
observations and numerical solutions) can be used to obtain the
estimates of the unknown parameters. The NLS method for (\ref
{Const_ODE}) was first described by
mathematicians in 1970s [\citet{r4}, \citet{r67}, \citet{r5}].
The NLS method was also widely used to estimate the
unknown parameters in ODE models in the fields of mathematics, computer
science and control engineering. In the 1990s, the NLS method was
extended to estimate time-varying parameters in (\ref
{Time_ODE}). For example, the NLS method with spline approximation to
time-varying parameters has been successfully applied to
pharmacokinetic [\citet{r33}], physiologic [\citet{r65}]
and HIV studies [\citet{r1}].

Though the NLS was the earliest and the most popular method developed
for estimating the parameters
in ODE models, so far the proposed NLS estimators and their asymptotic
properties for ODE models
have not been systematically studied, in particular, for time-varying
parameter estimates. The
influence of the numerical approximation error of ODEs on the
asymptotic properties has not been
analyzed. All existing studies took the numerical solution as the true
solution and did not consider
the difference between them. The difficulty is due to the co-existence
of both measurement error and
numerical error, and the standard theories of the NLS method [\citet{r32}, \citet{r37}, \citet{r72}, \citet{r14}] cannot be directly
applied. In this article, we intend to fill this gap.

The rest of the article is organized as follows. In Section \ref
{sec2}, we
discuss the identifiability problem of ODE models. Then we introduce
the numerical solution-based NLS estimators for (\ref
{Const_ODE}) and (\ref{Time_ODE}), and study their asymptotic
properties in Sections \ref{sec3} and \ref{sec4}, respectively. The
asymptotic properties
of the proposed estimators, including strong consistency, rate of
convergence and asymptotic normalities, are established using the tools
of empirical processes [Pollard (\citeyear{r47}, \citeyear{r49}), \citet{r43},
\citet{r66}, \citet{r36}, \citet{r71}]
and the sieve methods [\citet{r17}, \citet{r59}, \citet{r23} and
\citet{r58}]. We perform simulation studies to
investigate the finite-sample performance of the proposed estimation
methods in Section \ref{sec5}. In this section, we also apply the proposed
approaches to a set of ODE models for HIV dynamics. We provide a
summary and discussion for the proposed methods in Section \ref{sec6}. Finally,
the proofs of all the theoretical results are given in the \hyperref
[app]{Appendix}.

\section{Identifiability of ODE models}\label{sec2}

Identifiability of ODE models is a critical question to answer
before parameter estimation. To verify the uniqueness of parameter
estimates for given system inputs and outputs, both analytical and
numerical techniques have been developed for ODE models since 1950s.
Before jumping into technical details, two commonly used definitions
of identifiability are given as follows
[\citet{Bellman1970}, Cobelli, Lepschy and Jacur (\citeyear{Cobelli1979}), \citet{Walter1987},
\citet{Ljung1994}, \citet{Audoly2001}, \citet{r31}].
\begin{definition}
Globally identifiable: a system
structure is said to be globally identifiable if for any two
parameter vectors $\bolds{\beta}_1$ and $\bolds{\beta}_2$ in the
parameter space $\mathcal{B}$, $\mathbf{X}(t,\bolds{\beta}_1) =
\mathbf{X}(t,\bolds{\beta}_2)$ can be satisfied for all $t$ if and only
if $\bolds{\beta}_1 = \bolds{\beta}_2$.
\end{definition}

However, global identifiability is a strong condition to satisfy and
usually difficult to verify in practice. Therefore, the definition of
at-a-point identifiability was introduced by \citet{Ljung1994} and
\citet{Quaiser2009} as follows.
%
%
\begin{definition}
At-a-point identifiable: a system
structure is said to be locally (or globally) identifiable at a
point $\bolds{\beta}_*$ if for any $\bolds{\beta}$ within an open
neighborhood of $\bolds{\beta}_*$ (or within the entire parameter
space), $\mathbf{X}(t,\bolds{\beta})
= \mathbf{X}(t,\bolds{\beta}_*)$ can be satisfied for all $t$ if and only
if $\bolds{\beta} = \bolds{\beta}_*$.
\end{definition}

A number of methods have been proposed for identifiability analysis of
ODE models, including
\textit{structural} [\citet{Bellman1970}, \citet{Ljung1994}, \citet{XiaMoog2003}],
\textit{practical} [e.g., \citet{r54}, \citet{r39}] and
\textit{sensitivity-based} [e.g., \citet{Jolliffe1972}, \citet{Quaiser2009}] approaches. Due
to the limited space, we may not be able to provide an exhaustive list
of publications on
identifiability of ODE models. In this article, the structural
identifiability analysis techniques
are of particular interest mainly due to the theoretical completeness.

Various structural identifiability approaches have been proposed,
such as power series expansion [\citet{Pohjanpalo1978}], similarity
transformation [\citet{Vajada1989}, \citet{Chappel1992}] and
implicit function theorem method [\citet{Xia2003}, \citet{XiaMoog2003}, \citet{r39},
\citet{Wu2008}].
Particularly, \citet{Ollivier1990} and \citet{Ljung1994}
introduced another approach in the framework of differential algebra
[\citet{Ritt1950}, \citet{Kolchin1973}]. The differential algebra approach is
suitable to general nonlinear dynamic systems, and it has been
successfully applied to nonlinear differential equation models,
including models with time-varying parameters [\citet{Audoly2001}].
In this article, the differential algebra approach is employed to
verify the identifiability of ODE models with both constant and
time-varying parameters.

For most structural identifiability analysis techniques such as the
implicit function theorem method
and the differential algebra approach, a key step is the elimination of
latent variables via taking
derivatives and algebraic operations, which makes such techniques
suitable for multivariate ODE
models with partially observed state variables. After all unobserved
state variables are eliminated,
equations involving only given inputs, measured outputs and unknown
parameters can be obtained. If
we consider the parameters as unknowns, it is easy to verify that the
identifiability of unknown
parameters is determined by the number of roots of these equations.

For illustration purposes, we consider a classical HIV dynamic model
with both
constant and time-varying parameters [\citet{r42}, \citet{HRW03},
\citet{r74}] as an example:
%
%
\begin{equation} \label{HIV_model}
\cases{
\dfrac{d}{dt}T_U(t) =\lambda- \rho T_U(t)-\eta(t) T_U(t) V(t), \vspace*{2pt}\cr
\dfrac{d}{dt}T_I(t) =\eta(t) T_U (t) V(t)-\delta T_I(t), \vspace*{2pt}\cr
\dfrac{d}{dt}V(t) =N \delta T_I(t)-c V(t),}
\end{equation}
where $T_U$ is the concentration of uninfected target CD4$+$ $T$ cells,
$T_I$ the concentration of
infected CD4$+$ $T$ cells, $V(t)$ the viral load, $\lambda$ the
proliferation rate of uninfected CD4$+$
$T$ cells, $\rho$ the death rate of uninfected CD4$+$ $T$ cells, $\eta
(t)$ the time-varying infection
rate depending on antiviral drug efficacy, $\delta$ the death rate of
infected cells, $c$ the
clearance rate of free virions, $N$ the number of virions produced by a
single infected cell on
average. This model will also be used in our numerical studies in
Section \ref{sec5}. For notational
simplicity, let $x_1$, $x_2$ and $x_3$ denote $T_U$, $T_I$ and $V$, and
let $y_1=T_U+T_I=x_1+x_2$
and $y_2=V=x_3$ denote the measurable outputs, respectively. Then
(\ref{HIV_model}) can be re-written as
%
%
\begin{equation}\label{odex}
\cases{
x_1' =\lambda-\rho x_1-\eta(t) x_1 x_3, \vspace*{1pt}\cr
x_2' =\eta(t) x_1 x_3 - \delta x_2, \vspace*{1pt}\cr
x_3' =N \delta x_2-c x_3,}
\end{equation}
where $x_1', x_2'$ and $x_3'$ denote the derivatives of $x_1, x_2$ and
$x_3$, respectively. We adopt
the following ranking for variable elimination [\citet{Ljung1994}],
%
%
\begin{equation} \label{3DmodelRANKING}
\eta\prec y_2 \prec y_1 \prec\bolds{\beta} \prec x_3 \prec x_2
\prec x_1,
\end{equation}
where $\bolds{\beta} = (\lambda,\rho,N,\delta,c)^T$ is the vector of
constant unknown parameters. By taking the higher order derivatives
on both sides of (\ref{odex}) and using some algebra
elimination techniques, we can eliminate $x_1$, $x_2$ and $x_3$ from
(\ref{odex}) using the ranking (\ref{3DmodelRANKING}) to obtain
%
%
\begin{eqnarray}\label{E:CSetEq1} 
y_1^{(2)} + (\rho+ \delta) y_1' + \delta\rho y_1 -
\delta\lambda+ \eta(t) y_2 (y_1' + \delta y_1 - \lambda) = 0, \\
\label{E:CSetEq2}
y_2^{(2)} + (\delta+ c) y_2' + \delta c y_2 - \eta(t) y_2 (N \delta
y_1 -
y_2' - c y_2) = 0,
\end{eqnarray}
where $y_1^{(2)}$ and $ y_2^{(2)}$ denote the second-order derivative
of $y_1(t)$ and $y_2(t)$, respectively.
Therefore, $\eta(t)$ can be expressed in terms of measurable system
outputs and other constant unknown parameters
either from (\ref{E:CSetEq1}) as
%
%
\begin{equation}\label{ETA1}
\eta(t) = \frac{y_1^{(2)} + (\rho+ \delta) y_1' +
\delta\rho y_1 - \delta\lambda}{-y_2 (y_1' + \delta y_1 - \lambda)}
\end{equation}
or from (\ref{E:CSetEq2}) as
%
%
\begin{equation}\label{ETA2}
\eta(t) = \frac{y_2^{(2)} + (\delta+ c) y_2' +
\delta c y_2}{y_2 (N \delta y_1 - y_2' - c y_2)}.
\end{equation}
Thus, $\eta(t)$ is identifiable if all the constant unknown parameters
are identifiable. To verify
the identifiability of all unknown parameters $\bolds\theta=(\bolds
{\beta
}^T,\eta)^T$, equations (\ref{ETA1})
and (\ref{ETA2}) can be combined to obtain
%
%
\begin{eqnarray} \label{CUPID}\quad
&&y_1^{(2)} y_2 y_2' - y_1' y_2 y_2^{(2)} - \delta y_1 y_2 y_2^{(2)} +
\lambda y_2 y_2^{(2)} - (\delta+ c) y_1' y_2 y_2'\nonumber\\
&&\qquad{} +
(\rho\delta+ \rho+ \delta- \delta^2 - \delta c) y_1 y_2 y_2'
+ c y_2 y_2'\nonumber\\[-8pt]\\[-8pt]
&&\qquad{} + \rho c y_1' {y_2}^2 +
(\rho\delta c - \delta^2 c) y_1 {y_2}^2- N \delta y_1 y_1^{(2)}
y_2\nonumber\\
&&\qquad{} +c y_1^{(2)} {y_2}^2 - N \delta(\rho+ \delta) y_1 y_1' y_2
- N \delta^2 \rho{y_1}^2 y_2 + N \delta^2 \lambda y_1 y_2 = 0.\nonumber
\end{eqnarray}
The equation above only involves measurable system outputs
[$(T_U+T_I)$, $V$ and their
derivatives] and constant unknown parameters. We assume that the
derivatives of $(T_U+T_I)$ and $V$
exist and are continuous up to order 2. Although the derivatives of
$(T_U+T_I)$ and $V$ are usually
not directly measured in experiments, for theoretical identifiability
analysis, they are known once
$(T_U+T_I)$ and $V$ are measured (e.g., via numerical evaluation).
Finally, it can be verified that
(\ref{CUPID}) is of order 0 and of degree $>1$ in $\bolds
{\theta
}$, so (\ref{CUPID}) satisfies
the sufficient conditions given in \citet{Ljung1994} and $\bolds
{\beta
} = (\lambda,\rho,N,\delta,
c)^T$ is thus at-a-point identifiable at the true parameter point.
Therefore, $\eta(t)$ is also
at-a-point identifiable at the true parameter point. For more detailed
techniques for
identifiability analysis of ODE models, we refer readers to the
references listed above.

\section{ODE models with constant parameters}\label{sec3}
Throughout this article, we let $\|\mathbf{a}\|$ be the Euclidean norm
(or $L_2$ norm) of a vector (or a matrix) $\mathbf{a}$;
$\|\mathbf{A}\|_{\infty}={\max_{1\leq i\leq
m}\sum_{j=1}^n}|a_{ij}|$ be the supremum norm of an $m\times
n$ matrix $\mathbf{A}$, where $a_{ij}$ is the $(i,j)$th element of
$\mathbf{A}$; $\mathbf{A}^{\otimes2}=\mathbf{A}\mathbf{A}^T$ for a matrix
$\mathbf{A}$; $C^r[a,b]$ be the class of functions with $r$-order
continuous derivative on the interval $[a,b]$;
$\|f\|_{\infty}=\sup_t|f(t)|$ be the supremum norm of a function
$f$; and $x\wedge y$ denotes $\min(x,y)$. Moreover, for a random
vector $\mathbf{Z}\sim P$, where $P$ is a probability measure, we let
$\|f(\mathbf{Z})\|_2=\|f\|_{P,2}=(\int f^2 \,dP)^{1/2}$ be the
$L_2(P)$-norm of a function $f$.

\subsection{Measurement model and estimator}\label{sec31}

In this section, we consider ODE models with constant parameters,
that is, equation (\ref{Const_ODE}), over the time range of interest
$I=[t_0,T]$ ($-\infty<t_0<T<+\infty$), where the initial value
${\mathbf X}_0={\mathbf X}(t_0)$ is assumed to be known in this article.
In reality, ${\mathbf X}(t)$ cannot be measured exactly and directly;
instead, its surrogate $\mathbf{Y}(t)$ can be measured. For
simplicity, here we assume an additive measurement error model to
describe the relationship between ${\mathbf X}(t_i)$ and the surrogate
$\mathbf{Y}(t_i)$,
%
%
\begin{equation}\label{meas_err}
\mathbf{Y}(t_i)={\mathbf X}(t_i)+{\bolds\varepsilon}(t_i),
\end{equation}
at random or fixed design time points $t_1,\ldots,t_n$, where
the measurement errors $({\bolds\varepsilon}(t_1),\ldots,{\bolds
\varepsilon}(t_n))$ are independent with mean zero and a diagonal
variance--covariance matrix $\Sigma$. Moreover, in the case of random
design, assume that the measurement errors are independent of ${\mathbf
X}(t)$. Equation (\ref{meas_err}) is called the observation or
measurement equation.

If (\ref{Const_ODE}) does not have a closed-form solution, we
need to resort to numerical techniques to obtain numerical solutions
at discrete time points. In this article, we consider a general
one-step numerical method. Let $t_0=s_0<s_1<\cdots<s_{m-1}=T$ be
grid points on the interval $I$, $h_j=s_j-s_{j-1}$ be the step
size and $h=\max_{1\leq j\leq m-1}h_j$ be the maximum step
size, and\vspace*{1pt} ${\mathbf X}^h_j$ and ${\mathbf X}^h_{j+1}$ be the numerical
approximations to the true solutions ${\mathbf X}(s_j)$ and ${\mathbf
X}(s_{j+1})$, respectively, which can be typically written as
%
%
\begin{equation}\label{num_meth}
{\mathbf X}^h_{j+1}={\mathbf X}^h_j+h{\bolds\Phi}(s_j,{\mathbf
X}^h_j,{\mathbf X}^h_{j+1},h),
\end{equation}
where the specific form of $\bolds\Phi$ depends on the numerical method.
The common numerical methods
include the Euler backward method, the trapezoidal rule, the $r$-stage
Runge--Kutta algorithm ($r$ is
usually between 2 and 5), and so on. Among these algorithms, the
4-stage Runge--Kutta algorithm
[\citet{r38}, page 53, \citet{r18}, page 134] has been well
developed and widely used in practice. Therefore, we employ the 4-stage
Runge--Kutta algorithm as an
example in our numerical studies.

Define $e^h={\max_{0\leq j\leq m-1}}\|\mathbf{X}(s_j)-\mathbf{X}_j^h\|
$, which
is called the numerical error or the global discretization error
[\citet{r18}, page 159, \citet{r38}, page 57]. If
$e^h=O(h^p)$, $p$ is called the order of the numerical method. It is
necessary to establish a relationship between the number of grid
points $m$ (or the maximum step size $h$) and the sample size of
measurements $n$ since the asymptotic properties of the proposed
estimators are related to both numerical error and measurement
error. To our best knowledge, this is the first attempt to establish
such as a relationship.

Following Mattheij and Molenaar [(\citeyear{r38}), page 58] the interpolation
technique is commonly used if the measurement points
($t_i,i=1,2,\ldots,n$) are not coincident with the grid points
($s_j,j=1,2,\ldots,m-1$) of the numerical method, and the cubic Hermite
interpolation is often adopted. Let $\tilde{\mathbf X}(t,\bolds{\beta})$
denote the interpolated numerical solution of ${\mathbf X}(t,\bolds
{\beta})$ obtained from the numerical method for given
$\bolds{\beta}$, and then (\ref{meas_err}) can be approximately
rewritten as $\mathbf{Y}(t)\approx\tilde{\mathbf X}(t,\bolds{\beta
}_0)+{\bolds\varepsilon}(t)$. The simple numerical
solution-based NLS estimator $\hat{\bolds{\beta}}_n$ of
$\bolds{\beta}_0$ minimizes
%
%
\begin{equation}\label{obj_fun}
\Xi_1(\bolds\beta)=\sum_{i=1}^n\sum_{j=1}^K[Y_j(t_i)-\tilde
{\mathbf X}_j(t_i,\bolds{\beta})]^2.
\end{equation}
Note that if the data are correlated or the measurement variances
are heterogeneous, the weighted NLS can be used. The theoretical
results can be extended to the weighted NLS. Also note that we can
easily obtain the estimator $\hat{\mathbf X}(t)=\tilde{\mathbf
X}(t,\hat{\bolds{\beta}}_n)$ for ${\mathbf X}(t)$.

To minimize the NLS objective function (\ref{obj_fun}), the standard
gradient optimization methods
may fail due to the complicated nonlinear ODE model and the NLS
objective function may have multiple
local minima or may be ill-behaved [\citet{r16}].
Fortunately, various global
optimization methods are available to more reliably solve the parameter
estimation problem for ODE
models, although the global optimization methods are very
computationally intensive. \citet{r41} compared the performance and computational cost of
seven global optimization
methods, including the differential evolution method [\citet{r62}]. Their results suggest
that the differential evolution method outperforms the other six
methods with a reasonable
computational cost. Improved performance can be achieved using a hybrid
method combining gradient
methods and global optimization methods. A hybrid method based on the
scatter search and sequential
quadratic programming (SQP) has been proposed by \citet{r54}, who
showed that the hybrid scatter search method is much faster than the
differential evolution method
for a simple HIV ODE model. In addition, \citet{r39} also
suggested that global optimization
methods should be used for general nonlinear ODE models. Here we
combine the differential evolution, the scatter search method and the
SQP local optimization technique to implement our NLS minimization.

\subsection{Asymptotic properties}\label{sec32}
In this section, we study the asymptotic properties of the proposed
numerical solution-based NLS estimator when both measurement error and
numerical error are considered. First we make the following assumptions:

\begin{enumerate}[A13.]
\item[A1.] $\bolds\beta\in\mathcal{B}$, where $\mathcal{B}$ is a compact
subset of $\mathcal{R}^d$ with a finite diameter
$R_{\bolds\beta}$.
\item[A2.] $\Omega=\{{\mathbf X}(t,\bolds{\beta})\dvtx t\in I,
\bolds{\beta}\in
\mathcal{B}\}$ is a closed and bounded convex subset of $\mathcal{R}^K$.
\item[A3.] There exist two constants $-\infty<\mathbf{c}_{{1}}<
\mathbf{c}_{{2}}<+\infty$ such that $c_1\leq
\mathbf{Y}(t)\leq c_2$ for all $t\in I$.
\item[A4.] All partial derivatives of ${\mathbf F}(t,\mathbf
{X},\bolds\beta)$
up to order $p$ with respect to $t$ and $\mathbf{X}$
exist and are continuous.
\item[A5.] The numerical method for solving ODEs is of order $p$.
\item[A6.] For any
$\bolds{\beta}\in\mathcal{B}$, $E_t[{\mathbf X}(t,\bolds{\beta
})-{\mathbf X}(t,\bolds{\beta}_0)]^2=0$ if and only if
$\bolds{\beta}=\bolds{\beta}_0$.
\item[A7.] The first and second partial derivatives, $\frac{\partial
{\mathbf X}(t,\bolds\beta)}{\partial\bolds{\beta}}$ and
$\frac{\partial^2 {\mathbf X}(t,\bolds\beta)}{\partial\bolds
{\beta}\,\partial\bolds{\beta}^T}$, exist
and are continuous and uniformly bounded for all $t\in I$ and
$\bolds{\beta}\in\mathcal{B}$.
\item[A8.] For the ODE numerical solution $\tilde{\mathbf X}(t,\bolds
{\beta})$,
the first and second partial derivatives,
$\frac{\partial{\tilde{\mathbf{X}}}(t,\bolds\beta)}{\partial
\bolds{\beta}}$
and $\frac{\partial^2 {\tilde{
\mathbf{X}}}(t,\bolds\beta)}{\partial\bolds{\beta}\,\partial\bolds
{\beta}^T}$, exist
and are continuous and uniformly bounded for all $t\in I$ and
$\bolds{\beta}\in\mathcal{B}$.
\item[A9.] Let $0<c_3<c_4<\infty$ be two constants. For random design
points, $t_1, \ldots, t_n$ are
i.i.d. The joint density function $\phi(t,\mathbf{y})$ of $(t,\mathbf{Y})$
satisfies $c_3\leq\phi(t,\mathbf{y})\leq c_4$ for all $(t,\mathbf
{y})\in
[t_0,T]\times[\mathbf{c}_{{1}},\mathbf{c}_{{2}}]$.
\item[A10.] The true parameter ${\bolds\beta}_0$ is an interior
point of
$\mathcal{B}$.
\item[A11.] ${\tilde{\bolds\beta}}$ is an interior point of
$\mathcal
{B}$, where
$\tilde{\bolds\beta}=\arg\min_{\bolds\beta\in\mathcal
{B}}E_0[\mathbf{Y}(t)-\tilde{\mathbf{X}}(t,\bolds\beta)]^T\times\break[\mathbf
{Y}(t)-\tilde{\mathbf{X}}(t,\bolds\beta)]$ and $E_0$ is the expectation
with respect to $P_{\bolds\beta_0}$, the joint probability
distribution of $(t,\mathbf{Y}(t))$ at true value
$\bolds\beta_0$.\vspace*{1pt}
\item[A12.] $\mathbf{V}_1=\{E_t(\frac{\partial
{\mathbf X}}{\partial{\bolds\beta_0}}\,\frac{\partial{\mathbf
X}}{\partial{\bolds\beta_0^T}})\}^{-1}E_t(\frac{\partial{\mathbf
X}}{\partial{\bolds\beta_0}}\Sigma\frac{\partial{\mathbf
X}}{\partial{\bolds\beta_0^T}})\{E_t(\frac{\partial{\mathbf
X}}{\partial
{\bolds\beta_0}}\,\frac{\partial
{\mathbf X}}{\partial{\bolds\beta_0^T}})\}^{-1}$ is positive
definite, where
$E_t[g(t)]$ is expectation
of function $g(t)$ with respect to $t$.
\item[A13.] $\tilde{\mathbf{V}}_1=\{E_t(\frac{\partial
\tilde{\mathbf{X}}}{\partial{\tilde{\bolds\beta}}}\frac{\partial
\tilde{\mathbf{X}}}{\partial{\tilde{\bolds\beta}^T}})\}
^{-1}E_0(\frac
{\partial
\tilde{\mathbf{X}}}{\partial{\tilde{\bolds\beta}}}[\mathbf
{Y}(t)-\tilde{
\mathbf{{X}}}(t,\tilde{\bolds\beta})]^{\otimes
2}\frac{\partial
\tilde{\mathbf{X}}}{\partial{\tilde{\bolds\beta}^T}})\{E_t(\frac
{\partial
\tilde{\mathbf{X}}}{\partial{\tilde{\bolds\beta}}}\frac{\partial
\tilde{\mathbf{X}}}{\partial{\tilde{\bolds\beta}^T}})\}^{-1}$ is positive
definite.
\end{enumerate}

Assumptions A1--A4 are general requirements for existence of numerical
solutions of ODE models.
Assumption A5 from Mattheij and Molenaar (\citeyear{r38}, pages 55 and 56) defines
the precision of the numerical
algorithm. For example, the Euler backward method, the trapezoidal
rule, the 4-stage and 5-stage
Runge--Kutta are of order 1, 2, 4 and 5, respectively. Theorem 2.13 in
Hairer, N{\o}rsett and Wanner [(\citeyear{r18}), page 153] provides sufficient and necessary conditions
for the numerical method to be of
order $p$. Theorems 3.1 and 3.4 in Hairer, N{\o}rsett and Wanner [(\citeyear{r18}), pages 157 and 160] give the
magnitude of the numerical error of the numerical algorithms.
Assumption A6 is required for
identifiability and imposed for consistency. From Section \ref{sec2},
we know
that the HIV model
(\ref{HIV_model}) is at-a-point identifiable at the true value $\bolds
{\beta}_0$. This result and
assumption A9 are sufficient conditions for assumption A6 to be
satisfied. Assumptions A7--A9 are
needed for consistency. Assumptions A10--A13 are needed for the proof of
asymptotic normality in
Theorem \ref{Theorem32}.
\begin{theorem}\label{Theorem31} Assume that there exists a $\lambda>0$ such
that $h=O(n^{-\lambda})$, then under assumptions \textup{A1--A10}, we have
$\hat{\bolds\beta}_n-\bolds\beta_0\rightarrow0$, almost surely under
$P_{\bolds\beta_0}$.
\end{theorem}
\begin{theorem}\label{Theorem32}
\textup{(i)} For $h=O(n^{-\lambda})$ with $\lambda>1/(p\wedge4) $ where $p$
is the order of the numerical method (\ref{num_meth}), under
assumptions \textup{A1--A10} and \textup{A12}, we have that
$n^{1/2}(\hat{\bolds\beta}_n-{\bolds\beta}_0)\stackrel
{d}{\rightarrow}
N(0,\mathbf{V}_1)$.

\textup{(ii)} For $h=O(n^{-\lambda})$ with $0<\lambda\leq1/(p\wedge4)$,
under assumptions \textup{A1--A9}, \textup{A11} and \textup{A13}, we have that
$n^{1/2}(\hat{\bolds\beta}_n-{\tilde{\bolds\beta}})\stackrel
{d}{\rightarrow}
N(0,\tilde{\mathbf{V}}_1)$ with
$\|\tilde{\bolds\beta}-\bolds\beta_0\|=O(h^{(p\wedge
4)/2})=O(n^{-\lambda(p\wedge4)/2})$ and
$\|\tilde{\mathbf{V}}_1-V_1\|=O(h^{(p\wedge4)/2})=O(n^{-\lambda
(p\wedge4)/2})$.
\end{theorem}

The detailed proofs of Theorems \ref{Theorem31} and \ref{Theorem32} are
provided in the \hyperref[app]{Appendix}. The basic idea for the proofs
is motivated by \citet{r43} in which a general central limit
theorem is proved for a broad class of simulation estimators, that is,
the objective function of the estimator is too complicated to evaluate
directly, and instead the Monte Carlo simulation is used to approximate
the objective function to obtain the estimator. The asymptotic
properties of the simulation-based estimator are established using a
general central limit theorem under nonstandard conditions given in
Huber (\citeyear{r29}) and \citet{r48}, which are often called the
Huber--Pollard Z-theorem [see Theorem 3.3.1 in \citet{r66}]. In
this article, we use the same theorem to prove the asymptotic normality
of the numerical solution-based NLS estimator for ODEs. Similarly, our
objective function $\Xi_1(\beta)$ in (\ref {obj_fun}) cannot be
directly evaluated; instead we have to approximate it by solving
(\ref{Const_ODE}) numerically. Thus, similar ideas in
\citet{r43} can be borrowed to establish the asymptotic results of
our estimator in Theorems \ref{Theorem31} and \ref{Theorem32}.
\begin{Remark}\label{Remark1}
Theorems \ref{Theorem31} and \ref{Theorem32} can be extended to fixed
design points $t_i\in[t_0,T]$ ($i=1,\ldots,n$). Assume that there
exists a distribution function $Q(t)$ with corresponding density
$\varphi(t)$ such that, with $Q_n(t)$, the empirical distribution of
$(t_1,\ldots,t_n)$, ${\sup_{t\in [t_0,T]}}|Q_n(t)-Q(t)|=O_p(n^{-1/2})$
and $\varphi(t)$ is bounded away from zero and has continuous
second-order derivative on $[t_0,T]$. Define $E_t[g(t)]$ be the
integral $\int_{t_0}^T g(t)\,dQ(t)$ for function $g(t)$. Similarly we can
prove Theorems~\ref{Theorem31} and \ref{Theorem32} for the fixed design
if we replace assumption A9 by above assumption.
\end{Remark}
\begin{Remark}\label{Remark2}
From the proof of Theorem \ref{Theorem32} in the \hyperref[app]{Appendix}, we
still have
$\|\tilde{\bolds\beta}-\bolds\beta_0\|=O(h^{(p\wedge4)/2})$ and
$\|\tilde
{V}_1-V_1\|=O(h^{(p\wedge
4)/2})$ for a fixed constant $h$, which is independent of the sample
size $n$. This suggests that,
if the maximum step size $h$ of the numerical algorithm for solving
ODEs is a fixed constant, the
numerical solution-based NLS estimator is not consistent. Instead the
asymptotic bias is in the
order of $h^{(p\wedge4)/2}$.
\end{Remark}

Notice that our asymptotic results provide a
theoretical foundation for the relationship between the numerical
step size and sample size, that control numerical error and
measurement error, respectively, for the widely-used NLS estimator
based on the numerical solutions of the ODEs. Intuitively, the
smaller the numerical step size is, better the estimator is.
However, a smaller step size will increase the computational cost
and this may become a serious problem when the ODE system is large
and the computational cost is high. It is important to study the
trade-off between the numerical error and measurement error when the
computational cost needs to be taken into consideration. Our
theoretical results show that, only when the numerical step size,
which controls the numerical error and computational cost, goes to
zero with a rate faster than a particular rate $n^{-1/(p\wedge4)}$,
the numerical solution-based NLS estimator converges to the true
value of the parameters with the rate of root-$n$. In addition, the
asymptotic variance of the NLS estimator is the one as if the true
solution ${\mathbf X}(t)$ is exactly known.

The asymptotic variance--covariance matrix needs to be estimated in
order to perform statistical inference for unknown parameters
${\bolds\beta}$. There are some standard methods that can be used. The
first approach is to use the observed pseudo-information matrix based
on the NLS objective function (\ref{obj_fun}). The observed
pseudo-information matrix is defined as $\mathcal
{I}_1(\bolds\beta)=-\frac{\partial^2\Xi_1}{\partial\bolds\beta ^2}$.
The standard\vspace*{-2pt} error of $\hat{\bolds{\beta}}_{{n}}$
can then be approximated by
$\mathcal{I}_1^{-1/2}(\hat{\bolds\beta}_{{n}})/\sqrt{n}$. In
practice, we have noted that the inverse of the observed
pseudo-information matrix provides a reasonable approximation to the
asymptotic variance--covariance matrix $V_1$. \citet{r54} also
proposed this approach for parameter inference in ODE models.

The second approach is the weighted bootstrap method [\citet{r36}].
Let $W_i$, $i=1,\ldots,n$, denote $n$ i.i.d. positive
random weights with mean one [$E(W)=1$] and variance one
[$\operatorname{Var}(W)=1$]. The weights, $W_i$ are independent of
$\{\bolds\beta,t,\mathbf{Y}(t)\}$. For (\ref{Const_ODE}), the
weighted M-estimator $\hat{\bolds{\beta}}{}^0_n$ satisfies
\[
\hat{\bolds{\beta}}{}^0_n=\arg\min\sum_{i=1}^n\sum_{j=1}^K
W_i[Y_j(t_i)-\tilde{\mathbf{X}}_j(t_i,\bolds\beta)]^2.
\]
From Corollary 2 and Theorem 2 in \citet{r36}, given
$\{t_i,\mathbf{Y}(t_i)\}$,
$\sqrt{n}(\hat{\bolds{\beta}}{}^0_n-\hat{\bolds{\beta}}_n)$ and
$\sqrt{n}(\hat{\bolds{\beta}}_n-\bolds{\beta}_0)$ have the same limiting
distribution, then the weighted M-estimator $\hat{\bolds{\beta}}{}^0_n$
can be used for inference on $\hat{\bolds\beta}_n$.

Note that the empirical bootstrap has been used for statistical
inference for ODE models [\citet{r31a}]. However, the asymptotic
properties of the empirical
bootstrap estimators are quite difficult to derive. This is why we
propose to use the weighted
bootstrap method instead of the empirical bootstrap approach.

\section{ODE models with both constant and time-varying
parameters}\label{sec4}

\subsection{Measurement model and estimator}\label{sec41}

In this section, we consider (\ref{Time_ODE}) with both constant
and time-varying parameters, where the initial value ${\mathbf
X}_0={\mathbf X}(t_0)$ is assumed to be known. Again, ${\mathbf X}(t)$ is
not observed directly in practice; instead, we observe its surrogate
$\mathbf{Y}(t)$ through (\ref{meas_err}).

Let $\mathcal{A}$ be the following class of functions,
%
%
\begin{equation}\label{spline_def}
\mathcal{A}=\bigl\{\eta\in C^{\mu}[t_0,T]\dvtx
\bigl|\eta^{(\mu)}(z_1)-\eta^{(\mu)}(z_2)\bigr|\leq L|z_1-z_2|^{\gamma}\bigr\},
\end{equation}
where $\mu$ is a nonnegative integer, $\gamma\in(0,1]$,
$\varrho=\mu+\gamma>0.5$, and $L$ an unknown constant. The
smoothness assumption is often used in nonparametric curve
estimation. Usually, either $\varrho=1$ (i.e., $\mu=0$ and
$\gamma=1$) or $\varrho=2$ (i.e., $\mu=1$ and $\gamma=1$) should be
satisfied in various situations. Denote
$\bolds\theta=(\bolds{\beta}^T,\eta)^T$. Then the parameter space is
denoted by $\Theta=\{\bolds\theta\dvtx\bolds\beta\in\mathcal
{B},\eta\in\mathcal
{A}\}=\mathcal{B}\times\mathcal{A}$.

In this article, we use the method of sieves to approximate $\eta_0(t)$
on the support interval
$[t_0,T]$ of $t$. The basic idea of the sieve approach is to
approximate an infinite-dimensional
parameter space $\Theta$ by a series of finite-dimensional parameter
spaces $\Theta_n$, which depend
on the sample size $n$, and then to estimate the parameter on the
finite-dimensional spaces
$\Theta_n$ instead of $\Theta$. The concept of sieve was first proposed
by \citet{r17}. Since
then, the sieve method has been a powerful tool in the area of
nonparametric and semiparametric
statistics [\citet{r59}, \citet{r23}, \citet{r66}, Section 3.4,
\citet{r58}, \citet{r25}, \citet{r24}, \citet{r19},
\citet{r76} and \citet{r27}].

Here we apply the sieve estimation method to (\ref{Time_ODE})
with a time-varying parameter. First, we approximate $\eta(t)$ by
B-spline functions on the support interval $I$ of $t$. Let
$t_0=u_0<u_1<\cdots<u_q=T$ be a partition of the interval $I$, where
$q=O(n^v)$ ($0<v<0.5$) is a positive integer such that ${\max_{1\leq
j\leq q}}|u_j-u_{j-1}|=O(n^{-v})$. Then we have $N=q+l$ normalized
B-spline basis functions of order $l+1\geq\varrho$ [see \citet{r26}, page 1618]
that form a basis for the linear spline space. We
denote these basis functions in the forms of a vector
$\bolds{\pi}(t)=(B_1(t),\ldots,B_N(t))^T$ with which $\eta(t)$ can be
approximated by $\bolds{\pi}(t)^T{\bolds\alpha}$, where
${\bolds\alpha}=(\alpha_1,\ldots,\alpha_N)^T\in\mathcal{R}^N$ is the
spline coefficient vector with ${\bolds\alpha}_0$ corresponding to
$\eta_0(t)$. Such approximation is extensively used in nonparametric
and semiparametric problems [\citet{r61}, \citet{r59},
\citet{r58}, \citet{r24} and \citet{r26}]. The readers are referred to
Schumaker [(\citeyear{r57}), page 118] for more details about the construction of
the basis functions. Regression spline approximation to a
nonparametric function can always be expressed as a linear function
of basis functions so that the problem of time-varying coefficients
can be transformed into an estimation problem for a number of
constant parameters. Thus the estimation methods and computational
algorithms developed for (\ref{Const_ODE}) with constant
coefficients in Section~\ref{sec3} can be employed for (\ref
{Time_ODE})
with both constant and time-varying parameters.

For any $\bolds{\theta}_i\in\mathcal{B}\times\mathcal{A}$
($i=1,2$), we
define a distance
%
%
\begin{equation}\label{norm}
d(\bolds{\theta}_1,\bolds{\theta}_2)=\|\bolds{\beta}_1-\bolds
{\beta}_2\|+\|\eta
_1-\eta_2\|_2.
\end{equation}
Denote set
\[
\mathcal{A}_n=\Biggl\{\eta(t)=\sum_{i=1}^N B_i(t)\alpha_i\dvtx
\max_{1\leq i\leq N}|\alpha_i|\leq\ell_n\Biggr\},
\]
where $\ell_n\leq n^{(2l-1)/[2l'(2l+1)]}$ with a constant $l'$
arbitrarily close to $l$ [see \citet{r58}, page 2560],
then $\Theta_n=\{\bolds\theta\dvtx\bolds\beta\in\mathcal{B},
\eta\in\mathcal{A}_n\}=\mathcal{B}\times\mathcal{A}_n$ can be
used as a sieve
of $\Theta$. In fact, for any
$\bolds\theta=(\bolds\beta^T,\eta)^T\in\Theta$, by Corollary
6.21 in
\citet{r57}, there exists $\eta_n\in\mathcal{A}_n$ such that
$\|\eta_n-\eta\|_{\infty}=O_p(n^{-v\varrho})$. Denote
$\bolds\theta_n=(\bolds\beta^T,\eta_n)^T\in\Theta_n$, then
$d(\bolds\theta,\bolds\theta_n)=O_p(n^{-v\varrho})$. Equation
(\ref{Time_ODE}) now becomes
\[
\frac{d {\mathbf X}(t)}{dt}\approx{\mathbf F}\{t,{\mathbf
X}(t),\bolds{\beta},\pi(t)^T\bolds{\alpha}\}.
\]
For this approximation model, let $\tilde{\mathbf X}(t,\bolds{\beta
},\pi(t)^T\bolds{\alpha})$ be the numerical
approximation of ${\mathbf X}(t,\bolds{\beta},\eta(t))$ that can be
obtained from the same numerical algorithm as described in Section
\ref{sec3}. Equation (\ref{meas_err}) can be approximated by $\mathbf
{Y}(t)\approx\tilde{\mathbf X}(t,\bolds{\beta},\pi(t)^T\bolds
{\alpha}_0)+{\bolds\varepsilon}(t)$. The
numerical solution-based sieve NLS estimator
$\hat{\bolds{\theta}}_n=(\hat{\bolds\beta}{}^T_n,\hat{\eta}_n)^T$ is
defined as
%
%
\begin{equation}\label{esti}
\hat{\bolds{\theta}}_n=\mathop{\arg\inf}_{\bolds\theta\in
\Theta_n}\Xi_2(\bolds{\theta})
=\mathop{\arg\inf}_{\bolds\theta\in\Theta_n}
\sum_{i=1}^n\sum_{j=1}^K[Y_j(t_i)-\tilde
X_j(t_i,\bolds{\beta},\eta(t))]^2.
\end{equation}
When we substitute the sieve NLS estimators $\hat{\bolds{\theta}}_n$
into the numerical approximation, we can obtain the estimator
$\hat{\mathbf X}(t)=\tilde{\mathbf X}(t,\hat{\bolds{\beta}}_n,\hat
\eta_n(t))$.

\subsection{Asymptotic properties}\label{sec42}

The empirical objective function for the sieve NLS method proposed in
Section \ref{sec41} is a second-order loss function which is not a likelihood
function. We cannot use the standard information calculation of the
maximum likelihood estimator (MLE) based on orthogonal projections in
semiparametric models [\citet{r6}], and the asymptotic
normality theory for semiparametric MLEs obtained by Huang (\citeyear{r23},
Theorem 6.1) does not apply to our case. Fortunately,
\citet{r36} and \citet{r71} extended the Huang's asymptotic normality results to more
general semiparametric M-estimators by using a so-called
pseudo-information calculation. We are able to employ these new
asymptotic results to asymptotic properties of the proposed sieve NLS
estimator, and the following additional assumptions are needed:

\begin{enumerate}[B4.]
\item[B1.] The true time-varying parameter $\eta_0(\cdot)\in
\mathcal{A}$, where $\mathcal{A}$ is denoted in (\ref{spline_def}).
\item[B2.] All partial derivatives of ${\mathbf F}$ up to order $p$
with respect
to $t, {\mathbf X}$, and $\eta$, respectively, exist and are continuous.
\item[B3.] For any $\bolds{\beta}\in\mathcal{B}$ and $\eta\in
\mathcal{A}$,
$E_t[{\mathbf X}(t,\bolds{\beta},\eta(t))-{\mathbf X}(t,\bolds
{\beta}_0,\eta_0(t))]^2=0$ if and only if
${\bolds\beta}={\bolds\beta}_0$ and $P\{t\dvtx\eta(t)=\eta_0(t)\}=1$.
\item[B4.] The first and second partial
Fr\'{e}chet-derivatives [\citet{r66}, page 373]
in the norm $d$ defined in (\ref{norm}),
$\frac{\partial{\mathbf X}(t,\bolds\beta,\eta)}{\partial\bolds
{\beta}}$,
$\frac{\partial{\mathbf X}(t,\bolds\beta,\eta)}{\partial\eta}$,
$\frac{\partial^2 {\mathbf X}(t,\bolds\beta,\eta)}{\partial\bolds
{\beta}\,\partial\bolds{\beta}^T}$,
$\frac{\partial^2 {\mathbf X}(t,\bolds\beta,\eta)}{\partial\bolds
{\beta}\,\partial\eta}$ and
$\frac{\partial^2 {\mathbf X}(t,\bolds\beta,\eta)}{\partial\eta
^2}$ exist
and are continuous and uniformly bounded for all $t\in I$,
$\bolds{\beta}\in\mathcal{B}$ and $\eta\in\mathcal{A}$.
\item[B5.] For the ODE numerical solution $\tilde{{\mathbf
X}}(t,\bolds\beta
,\eta)$, the first and second partial
Fr\'{e}chet-derivatives in the norm $d$,
$\frac{\partial\tilde{\mathbf X}(t,\bolds\beta,\eta)}{\partial
\bolds{\beta
}}$, $\frac{\partial\tilde{\mathbf X}(t,\bolds\beta,\eta
)}{\partial\eta}$,
$\frac{\partial^2 \tilde{\mathbf X}(t,\bolds\beta,\eta)}{\partial
\bolds{\beta}\,\partial\bolds{\beta}^T}$,
$\frac{\partial^2 \tilde{\mathbf X}(t,\bolds\beta,\eta)}{\partial
\bolds{\beta}\,\partial\eta}$ and
$\frac{\partial^2 \tilde{\mathbf X}(t,\bolds\beta,\eta)}{\partial
\eta^2}$
exist and are continuous and uniformly bounded for all $t\in I$,
$\bolds{\beta}\in\mathcal{B}$ and $\eta\in\mathcal{A}$.
\item[B6.] For $K\geq2$, $\mathbf{V}_2=\mathbf{S}_1^{-1}\mathbf{S}_2(
\mathbf{S}_1^{-1})^T$ is positive
definite, where $\mathbf{S}_1$ and $\mathbf{S}_2$ are defined in
(\ref{S1})
and (\ref{S2}) in the \hyperref[app]{Appendix},
respectively.
\item[B7.] $v$ satisfies the restrictions $0.25/\varrho<v<0.5$ and
$v(2+\varrho)>0.5$, where $\varrho$ is the measure
of smoothness of $\eta(t)$ defined in assumption (B1).
\end{enumerate}
\begin{theorem}\label{Theorem41} Assume that there exists a $\lambda>0$ such that
$h=O(n^{-\lambda})$ and under assumptions \textup{A1--A4, A9, A10} and \textup{B1--B5},
we have $d(\hat{\bolds\theta}_n,\bolds\theta_0)\rightarrow0$, almost
surely under $P_{\bolds{\theta}_0}$.
\end{theorem}
\begin{theorem}\label{Theorem42} Assume that there exists a
$\lambda>1/[2(p\wedge4)]$ such that $h=O(n^{-\lambda})$ where $p$ is
the order of the numerical
algorithm (\ref{num_meth}), and under assumptions \textup{A1--A4, A9, A10} and
\textup{B1--B5}, we have
$d(\hat{\bolds\theta}_n,\bolds\theta_0)=O_p(n^{-v\varrho}+n^{-(1-v)/2})$.
\end{theorem}

From Theorem \ref{Theorem42}, we know that
$\|\hat{\bolds\beta}_n-\bolds\beta_0\|=O_p(n^{-v\varrho}+n^{-(1-v)/2})$
and
$\|\hat{\eta}_n(t)-\eta_0(t)\|_2=O_p(n^{-v\varrho}+n^{-(1-v)/2})$.
If $v=1/(1+2\varrho)$, the rate of convergence of $\hat{\eta}_n$ is
$n^{-\varrho/(1+2\varrho)}$, which is the same as the optimal rate
of the standard nonparametric function estimation [\citet{r60}].
Theorem \ref{Theorem43} below states that the rate of weak convergence of
$\hat{\bolds\beta}_n$ achieves $n^{-1/2}$ under some additional
assumptions.
\begin{theorem}\label{Theorem43}
For the maximum step size $h=O(n^{-\lambda})$ with
$\lambda>1/(p\wedge4) $, under assumptions \textup{A1--A4, A9, A10} and
\textup{B1--B7}, and $K\geq2$, we have
$n^{1/2}(\hat{\bolds\beta}_n-{\bolds\beta}_0)\stackrel
{d}{\rightarrow}
N(0,\mathbf{V}_2)$.
\end{theorem}
\begin{Remark}\label{Remark3}
For the case $h=O(n^{-\lambda})$ with $1/[2(p\wedge
4)]<\lambda\leq1/(p\wedge4)$, similar results to case (ii) in
Theorem \ref{Theorem32} can be obtained.
\end{Remark}

For $K=1$, Theorem \ref{Theorem43} does not hold, since in this case the special
perturbation direction $a^*(t)$ given in (\ref{simple_at}) is
$\frac{\partial X}{\partial\bolds\beta_0}/\frac{\partial
X}{\partial\eta_0}$, which leads to both $\mathbf{S}_1$ in (\ref{S1})
and $\mathbf{S}_2$ in (\ref{S2}) to be zero (see the proof of Theorem
\ref{Theorem43} in the \hyperref[app]{Appendix}). In this article, we consider one
special case that
we assume there exists an additive relationship between
$\bolds{\beta}$ and $\eta(\cdot)$ as follows:
%
%
\begin{equation}\label{additive}
\frac{d X(t)}{d t}=F\{t,X(t),\beta+\eta(t)\},
\end{equation}
which is a special case of (\ref{Time_ODE}), then the function
$X(t)$ has the form of
$X(t,\beta+\eta(t))$. In this case, we are able to establish similar
asymptotic normality results
under the identifiability constraint $E_t\eta(t)=0$. Note that \citet{r56} studied a similar
problem under a semiparametric regression model and used the same
identifiability constraint for the
unknown function $\eta(t)$ to establish the asymptotic normality for
the constant parameters. We
follow a similar idea and use B-spline approximation for $\eta(t)$. We
center the B-spline estimator
of $\eta(t)$ as follows:
\[
\hat\eta_n(t_i)\approx\sum_{i=1}^N
B_j(t_i)\hat\alpha_j-\frac1n\sum_{i=1}^n\sum_{j=1}^N
B_j(t_i)\hat\alpha_j=\sum_{j=1}^N\hat\alpha_j\Biggl[B_j(t_i)-\frac
1n\sum_{i=1}^n
B_j(t_i)\Biggr],
\]
which is subject to the constraints $\sum_{i=1}^n
\hat\eta_n(t_i)=0$. Under similar assumptions, the strong
consistency and the rate of weak convergence of the estimators,
similar to those of Theorems \ref{Theorem41} and \ref{Theorem42}, can be obtained. In
particular, the asymptotic normality can be established as follows:
\begin{Proposition}\label{Proposition1}
For (\ref{additive}) with
$K=1$, when the maximum step size $h=O(n^{-\lambda})$ with
$\lambda>1/(p\wedge4) $, under assumptions\vspace*{-1pt} \textup{A1--A4, A9, A10, B1--B5,
B7} and in addition $E_t[\eta(t)]=0$, we have
$n^{1/2}(\hat{\bolds\beta}_n-{\bolds\beta}_0)\stackrel
{d}{\rightarrow}
N(0,V_3)$, where $V_3=\sigma_0^2\{E_t(\frac{\partial
X}{\partial{\bolds\xi}})^ 2\}^{-1}$ with $\xi=\beta_0+\eta_0(t)$.
\end{Proposition}

The proof of this proposition is different from that of Theorem \ref{Theorem43}
and is given in the \hyperref[app]{Appendix}.
\begin{Remark}\label{Remark4}
By combining Theorem \ref{Theorem43} and Proposition \ref{Proposition1}, we can see
that the proposed sieve NLS
estimator is asymptotically normal with a convergence rate of $\sqrt
{n}$ for $K\geq2$ under
assumption B6, but we are only able to prove the result for a special
ODE model
(\ref{additive}) for $K=1$. This is because the asymptotic covariance
$V_2$ defined in B6 is always
singular in the case of $K=1$, and is only possibly nonsingular in the
case of $K\geq2$. Since
$V_2$ is always singular for $K=1$, we derive the asymptotic
distribution for the special ODE model
(\ref{additive}) using a different approach which results in
Proposition \ref{Proposition1}.
\end{Remark}

Similar approaches proposed in Section \ref{sec3} can be used to
estimate the
asymptotic variance--covariance matrix for $(\hat{\bolds\beta}_n,
\hat\eta_n(t))$. For the first approach, the observed
pseudo-information matrix can be evaluated by replacing $\eta(t)$
with the spline approximation $\pi^T(t)\bolds{\alpha}$, that is, to
rewrite the objective function $\Xi_2(\bolds{\theta})$ in the
expression (\ref{esti}) as $\Xi_2(\bolds{\beta},\bolds{\alpha})$. Then
the observed pseudo-information matrix $\bolds{\mathcal{I}}_2(\bolds
{\beta
},\bolds{\alpha})$ can be defined as
\[
\bolds{\mathcal{I}}_2(\bolds\beta,\bolds\alpha) =\pmatrix{
-\dfrac{\partial^2\Xi_2}{\partial\bolds\beta^2}&-\dfrac{\partial
^2\Xi
_2}{\partial\bolds\beta\,
\partial\bolds\alpha}\vspace*{2pt}\cr
-\dfrac{\partial^2\Xi_2}{\partial\bolds\alpha\,\partial\bolds\beta
}&-\dfrac
{\partial^2\Xi_2}{\partial\bolds\alpha^2}}.
\]
The standard error of $(\hat{\bolds\beta}_{{n}},\hat{\bolds
\alpha}_{{n}})$ is
approximately $\bolds{\mathcal{I}}_2^{-1/2}(\hat{\bolds\beta
}_{{n}},\hat{\bolds\alpha}_{{n}})/\sqrt{n}$ from
which the standard error of $\hat{\bolds\beta}_{{n}}$ can be
obtained. We
also find that the inverse of the observed pseudo-information matrix
provides a reasonable approximation to $\mathbf{V}_2$ via our
simulation studies in the next section.

Similarly the weighted bootstrap method can also be used. For
(\ref{Time_ODE}), the weighted M-estimators
$(\hat{\bolds{\beta}}{}^0_n,\hat{\bolds\alpha}^0_n)$ satisfy
\[
(\hat{\bolds{\beta}}{}^0_n,\hat{\bolds{\alpha}}^0_n)=\arg\min\sum
_{i=1}^n\sum_{j=1}^K
W_i[Y_j(t_i)-
\tilde{\mathbf{X}}_j(t_i,\bolds{\beta},\bolds{\pi}(t)^T\bolds
{\alpha})]^2.
\]
Based on Corollary 2 and Theorem 2 in \citet{r36}, given
$\{t_i,\mathbf{Y}(t_i)\}$,
$\sqrt{n}(\hat{\bolds{\beta}}{}^0_n-\hat{\bolds{\beta}}_n)$ and
$\sqrt{n}(\hat{\bolds{\beta}}_n-\bolds{\beta}_0)$ have the same limiting
distribution which can be used to justify the weighted bootstrap for
inference on $\hat{\bolds{\beta}}_n$ and $\hat{\eta}_n(t)$.

\section{Numerical studies}\label{sec5}

In this section, we consider the HIV dynamic model described in Section
\ref{sec2}. Recall that in this
system, $T_U(t)$, $T_I(t)$ and $V(t)$ are state variables and $(\lambda
,\rho,\delta,N,c,\eta(t))^T$
are kinetic parameters. By introducing the time-varying infection rate
$\eta(t)$ in this HIV
dynamic model, the model can flexibly describe the long-term viral
dynamics. In clinical studies,
only viral load, $V(t)$ and total CD4$+$ $T$ cell count,
$T(t)=T_U(t)+T_I(t)$, are closely monitored and
measured over time. For easy illustration and computational simplicity,
we fix the parameters
$\rho$ and $\delta$ in our numerical studies, and our objective is to
estimate three constant
parameters and one time-varying parameter, $(\lambda,N,c,\eta(t))^T$
based on measurements of viral
load and total CD4$+$ $T$ cell count.

\subsection{Monte Carlo simulation study}\label{sec51}

The following parameter values and initial conditions were used to
simulate observation data for (\ref{HIV_model}): $T_U (0) =
600$, $T_I (0) = 30$, $V(0) = 10^5$, $\lambda= 36$, $\rho= 0.108$,
$N=1000$, $\delta= 0.5$, $c=3$. For comparison purpose, we
generated the measurement data of $V(t)$ and $T(t)$ for four
scenarios in our simulation studies: (i) $\eta(t)=\eta$ is a small
constant, $\eta= 9.5\times10^{-6}$; (ii)~$\eta(t)$ is
time-varying but with a smaller (10\%) variation, $\eta(t) = 9\times
10^{-5} \times\{1-0.9 \cos(\pi t/400) \}$; (iii)
$\eta(t)=\eta$ is a larger constant, $\eta= 3.84\times10^{-5}$;
and (iv) $\eta(t)$ is time-varying but with a large (10-fold)
variation, $\eta(t) = 9\times10^{-5} \times\{1-0.9 \cos(\pi
t/40) \}$. Note that for cases (i) and (iii), the values of
constant $\eta$ were chosen to be approximately the average of
$\eta(t)$ over the period of time interval for cases (ii) and (iv),
respectively.

Let $y_1=T=T_U + T_I$ denote the total number of infected and
uninfected CD4$+$ $T$ cells and $y_2=V$
denote the viral load, the measurement models are given as follows:
\begin{eqnarray*}
y_{1i} &=& T(t_i) + \varepsilon_{1i},\\
y_{2i} &=& V(t_i) + \varepsilon_{2i},
\end{eqnarray*}
where $\varepsilon_{1i}$ and $\varepsilon_{2i}$ are independent and
follow normal distributions with mean
zero and variances $\sigma_{1i}^2$ and $\sigma_{2i}^2$, respectively.
The HIV dynamic model was
numerically solved within the time range $[0,20]$ to generate the
simulated data at each time
interval of 0.5 using the 4-stage Runge--Kutta algorithm. Consequently,
the corresponding sample size
is 40. The 20\% measurement errors were added to the numerical results
of the ODE model according to
the observation equations above. We applied the proposed estimation
methods in Sections \ref{sec3} and \ref{sec4} to
the simulated data for the 4 cases to evaluate the performance of the
proposed estimators and the
effect of the model misspecification. To stabilize the computational
algorithm, we log-transformed
the data. We also fixed parameters $\rho$ and $\delta$ as their true values.

For evaluating the performance of the estimation methods, we define
the average relative estimation error (ARE) as
\[
\mathrm{ARE}=\frac{1}{M} \sum_{j=1}^{M} \frac{ |\hat{\theta}_j-\theta|} {
|\theta|}
\times100\%,
\]
where $\hat{\theta}_j$ is the estimate of the parameter vector
$\theta$
from the $j$th simulation
data set, and $M=500$ is the total number of simulation runs.

%
%
\begin{sidewaystable}
\tablewidth=\textheight
\tablewidth=\textwidth
\caption{Simulation results for constant $\eta$ and the time-varying
$\eta(t)$ models. The ARE is
calculated based on 500 simulation runs for the HIV dynamic model. In
addition, $\sigma^2_{\mathrm{ODE}}$ is
the average of the estimated variance by the observed
pseudo-information, and $\sigma^2_{\mathrm{emp}}$ is the
empirical variance based on simulations. The sample size is $n=40$ and
the noise level is about
20\%}
\label{AREComparison}
\begin{tabular*}{\tablewidth}{@{\extracolsep{\fill}}l@{}ccd{2.2}d{2.2}d{2.2}d{3.1}ccd{4.1}d{4.3}d{4.3}@{}}
\hline
\multirow{2}{35pt}[-8pt]{\centering{\hspace*{-4.15pt}\textbf{Change\break
\hspace*{-4.15pt}of}
$\bolds{\eta(t)}$}}
& \multirow{2}{40pt}[-8pt]{\centering{\textbf{True} $\bolds{\eta(t)}$
\textbf{model}}}
& \multirow{2}{45pt}[-8pt]{\centering{\textbf{Fitted} $\bolds{\eta(t)}$
\textbf{model}}}
& \multicolumn{3}{c}{$\bolds\lambda$} & \multicolumn{3}{c}{$\bolds N$}
& \multicolumn{3}{c@{}}{$\bolds c$} \\[-4pt]
&  &
& \multicolumn{3}{c}{\hrulefill} & \multicolumn{3}{c}{\hrulefill}
& \multicolumn{3}{r@{}}{\hrulefill} \\
&  & & \multicolumn{1}{c}{\textbf{ARE(\%)}}
& \multicolumn{1}{c}{$\bolds{\sigma^2_{\mathrm{ODE}}}$} & \multicolumn{1}{c}{$\bolds{\sigma
^2_{\mathrm{emp}}}$} & \multicolumn{1}{c}{\textbf{ARE(\%)}}
& \multicolumn{1}{c}{$\bolds{\sigma^2_{\mathrm{ODE}}}$}
& \multicolumn{1}{c}{$\bolds{\sigma^2_{\mathrm{emp}}}$}
& \multicolumn{1}{c}{\textbf{ARE(\%)}} & \multicolumn{1}{c}{$\bolds{\sigma^2_{\mathrm{ODE}}}$}
& \multicolumn{1}{c@{}}{$\bolds{\sigma^2_{\mathrm{emp}}}$} \\
\hline
Small &Constant & Constant & 3.19 & 2.49 & 1.91 & 17.7 & 3.23e$+$04 & 4.94e$+$04
& 17.4 & 0.313 & 0.425 \\
& & Time-varying & 6.45 & 9.82 & 8.77 & 22.9 & 7.14e$+$04 &
8.63e$+$04 & 20.5 & 0.593 & 0.635 \\
&Time-varying & Constant & 3.77 & 2.38 & 2.08 & 17.9 & 3.36e$+$04 &
4.71e$+$04 & 19.8 & 0.331 & 0.432 \\
& & Time-varying & 6.40 & 9.16 & 8.55 & 22.6 & 6.53e$+$04 & 7.93e$+$04 &
20.9 & 0.543 & 0.637 \\
[4pt]
Large &Constant & Constant & 6.29 & 8.19 & 12.1 & 72.5 & 1.13e$+$06 & 8.53e$+$05
& 67.3 & 9.22 & 6.75 \\
& & Time-varying & 7.34 & 9.19 & 15.6 & 88.8 & 3.25e$+$06 &
1.13e$+$06 & 82.5 & 26.2 & 8.96 \\
&Time-varying & Constant & 94.2 & 13.7 & 7.02 & 994 & 5.86e$+$07 &
1.25e$+$08 & 1899 & 1660 & 3780 \\
& & Time-varying & 15.6 & 31.5 & 48.2 & 29.5 & 1.67e$+$05 & 1.67e$+$05 &
25.1 & 1.81 & 1.37 \\
\hline
\end{tabular*}
\end{sidewaystable}

In Table \ref{AREComparison}, the AREs of the constant parameters $(\lambda,N,c)$ are
listed. In addition, we also
report $\sigma^2_{\mathrm{ODE}}$ as the average of the estimated variance by the
observed pseudo-information
matrix and $\sigma^2_{\mathrm{emp}}$ as the empirical variance based on
simulation runs. Based on these
results, we can see that, when the change of $\eta(t)$ is small as a
function of time $t$ or $\eta$
is a small constant, the estimation of parameters is always good by
fitting a constant $\eta$ model
as observed by the low ARE values. However, when the change of $\eta
(t)$ is large or $\eta$ is a
large constant, misspecification of $\eta(t)$ may produce large AREs
for all parameter estimates. In
particular, when $\eta(t)$ is time-varying with a large variation,
using a constant $\eta$ model may
result in very poor estimates for all constant parameters. The variance
estimates based on the
pseudo-information agree well with the empirical estimates based on
simulations, which shows that
the pseudo-information-based variance estimate is reasonably good. The
evaluation of the bootstrap
variance estimation is prohibited in our simulation study due to high
computational cost.

%
\begin{figure}

\includegraphics{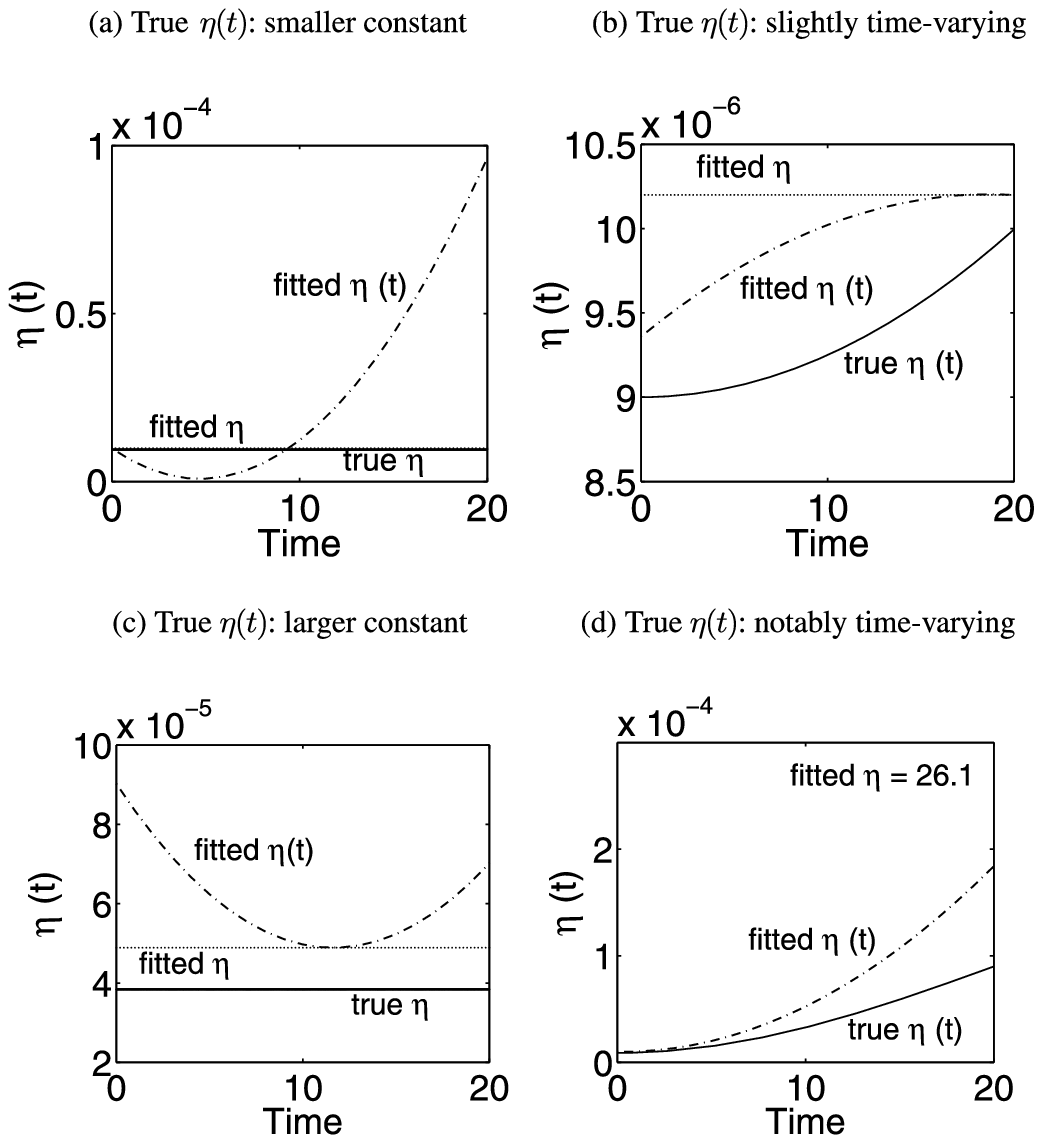}

\caption{Simulation
results for constant $\eta$ and the time-varying $\eta(t)$ models. In
each figure, the true model of
$\eta$ (solid), the constant $\eta$ model (dotted) and the time-varying
$\eta(t)$ model
(dash-dotted) are plotted and compared.} \label{figure1}
\end{figure}

In Figure \ref{figure1}, the average trajectories of estimated $\eta(t)$ are
compared to the true trajectories of
$\eta(t)$ for four different scenarios. From this figure, we observed a
similar trend as the
constant parameter estimates. The misspecification of $\eta(t)$
produces estimation error, in
particular for the cases with a large variation of $\eta(t)$ or a large
constant $\eta$. When the
model of $\eta(t)$ is correctly specified, the estimates based on the
proposed methods are
reasonably good. In order to evaluate the robustness of the proposed
approach, we also performed
further simulation studies for a complex function $\eta(t)=9.0\times
10^{-6}+9.0\times10^{-7}\times
t\times\{1-0.5\sin(\pi t/5.8)\}$ under the same simulation settings
(i.e., 40 time points, 20\%
error, 500 simulation runs). The results suggest that the sieve
estimator can still capture the
essential pattern of the complex $\eta(t)$ reasonably well (plots not shown).

\subsection{Application to AIDS clinical data}\label{sec52}

To illustrate applicability and feasibility of our proposed methods
and theories, we also applied the proposed estimation methods to fit
the HIV dynamic model to a clinical data set obtained from an HIV-1
infected patient who was treated with an antiretroviral therapy.
Very frequent viral load measurements were collected from this
patient after initiating the antiretroviral regimen: 13 measurements
during the first day, 14 measurements from day 2 to week 2, and then
one measurement at weeks 4, 8, 12, 14, 20, 24, 28, 32, 36, 40, 44,
48, 52, 56, 64, 74 and 76, respectively. In addition, the
measurements of total CD4$+$ $T$ cell counts were also taken at Day 1,
weeks 2 and 4, and monthly thereafter. Equation (\ref{HIV_model}) was
used to estimate HIV kinetic parameters using the viral load and
total CD4$+$ $T$ cell data.

For simplicity of illustration and computation, we fixed the initial
conditions of the state
variables in (\ref{HIV_model}) as $T_U (0) = 1$, $T_I (0) =
551$, $V(0) = 6.38\times10^4$,
which were derived from the baseline measurements. We also fixed two
parameters, as in the simulation
study, $\rho= 0.10$ and $\delta= 0.434$, which were taken from the
estimates in literature. Our
objective is to estimate the three constant parameters $(\lambda,N,c)$
and the time-varying
parameter $\eta(t)$ as in the simulation study. As we proposed in
Section \ref{sec4}, we employed B-splines
to approximate $\eta(t)$. We positioned the spline knots at
equally-spaced time points (the log-time
scale was used since the distribution of observation time points is
highly-skewed). We selected the
order of splines and the number of spline knots using the model
selection criterion AICc given by
\[
\mathrm{AICc} = n\ln\biggl(\frac{\mathrm{RSS}}{n}\biggr)+\frac{2nk}{n-k-1},
\]
where RSS is the residual of the sum of squares obtained from the NLS
model fitting, $n$ is the
total number of observations and $k$ is the number of unknown
parameters [including the
coefficients in the spline representation of $\eta(t)$]. Note that as a
practical guideline, if the
number of unknown parameters exceeds $n/40$ (where $n$ is the sample
size), the AICc instead of AIC
should be used. For our clinical data, the sample size $n$ is equal to
65, and the number of unknown
parameters varies between 6 and 13 for different scenarios, which is
much larger than $n/40 = 65/40
= 1.6$. Thus the AICc is more appropriate for our applications. In
general, the AICc converges to
the AIC as the sample size gets larger, thus the AICc is often
suggested to be employed regardless
of the sample size [\citet{r9}]. For our application,
we used AICc and compared the
models with the splines of order 3 and 4, and the number of knots from
3 to 10. In Table \ref{table2}, the AICc
values for these different models are reported, from which the best
model was selected as the spline
with order 3 and 5 knots for $\eta(t)$ approximation.

%
%
\begin{table}
\tablewidth=230pt
\caption{Model selection results for B-spline approximation of the
time-varying parameter $\eta(t)$} \label{table2}
\begin{tabular*}{\tablewidth}{@{\extracolsep{4in minus 4in}}lcd{2.0}c@{}}
\hline
\multicolumn{1}{@{}l}{\textbf{Model}} &
\multicolumn{1}{c}{\textbf{Spline}}
& \multicolumn{1}{c}{\textbf{Number of}} & \multicolumn{1}{c@{}}{\textbf{AICc}} \\
& \multicolumn{1}{c}{\textbf{order}} & \multicolumn{1}{c}{\textbf{knots}} &  \\
\hline
\phantom{0}1 & & 3 & $-$222.3 \\
\phantom{0}2 & & 4 & $-$242.8 \\
\phantom{0}3 & & 5 & $-$\textit{252.8} \\
\phantom{0}4 & & 6 & $-$243.8 \\
\phantom{0}5 & 3 & 7 & $-$250.6 \\
\phantom{0}6 & & 8 & $-$246.0 \\
\phantom{0}7 & & 9 & $-$246.8 \\
\phantom{0}8 & & 10 & $-$244.4 \\
[2pt]
\phantom{0}9 & & 3 & \multicolumn{1}{c@{}}{--} \\
10 & & 4 & $-$233.2 \\
11 & & 5 & $-$230.3 \\
12 & & 6 & $-$242.6 \\
13 & 4 & 7 & $-$249.1 \\
14 & & 8 & $-$245.5 \\
15 & & 9 & $-$244.9 \\
16 & & 10 & $-$240.5 \\
\hline
\end{tabular*}
\end{table}

We used the weighted bootstrap method to calculate both the confidence
intervals for the constant
parameters and the confidence bands for the time-varying parameter. The
basic idea of the weighted
bootstrap method is provided in Sections~\ref{sec32} and \ref{sec42}.
For the
computational implementation, we first
generated a positive random weight for each data point in the raw data
set from the exponential
distribution with mean one and variance one. By repeating this step, a
large number of (say, 1000)
sets of weights can be generated. Second, for each set of weights, the
ODE model is fitted to the
data to obtain parameter estimates by minimizing the weighted residual
sum of squares (see Sections
\ref{sec3} and \ref{sec4}). Recall that the time-varying parameter in
the model has been
approximated by B-splines,
then both the constant parameters and the constant B-spline
coefficients are actually estimated.
Once the estimates of the B-spline coefficients are obtained, we
construct the B-splines which
approximate the time-varying parameter. Thus, we eventually obtain 1000
estimates for each constant
parameter and 1000 B-splines for each time-varying parameter. Third,
for each constant parameter, we
select the 2.5\% and 97.5\% quantiles of the 1000 estimates to form the
95\% confidence intervals
for this parameter. For the time-varying parameter, at a single time
point, the 1000 B-splines have
1000 values. We also select the 2.5\% and 97.5\% quantiles of the 1000
values at this time point to
eventually form the 95\% pointwise confidence bands for the
time-varying parameter.

%
\begin{figure}

\includegraphics{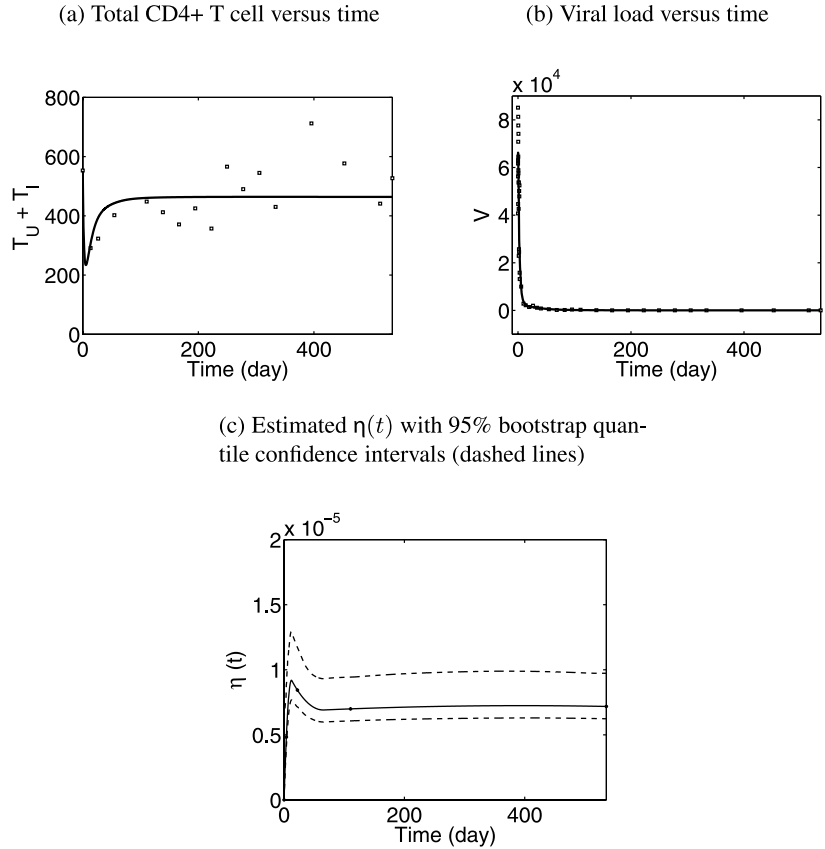}

\caption{Model fitting results with $\eta(t)$ approximated by
B-splines of order 3 and 5 knots.} \label{figure2}
\end{figure}

%
%
\begin{table}[b]
\tablewidth=250pt
\caption{The constant parameter estimation results}
\label{table3}
\begin{tabular*}{\tablewidth}{@{\extracolsep{\fill}}ld{4.2}c@{}}
\hline
\textbf{Parameter} & \multicolumn{1}{c}{\textbf{Estimate}} & \textbf{95\% confidence interval} \\
\hline
$\lambda$ & 46.52 & [43.20, 51.04] \\
$N$ & 1300.39 & [251.93, 4628.26] \\
$c$ & 4.35 & [0.98, 14.83] \\
\hline
\end{tabular*}
\end{table}

Model fitting results are given in Figure \ref{figure2} and Table \ref{table3}. From Figures
\ref{figure2}(a) and~(b), we can see that the
fitting is reasonably good for both CD4$+$ $T$ cell counts and viral load
data. The estimates of
constant parameters $(\lambda,N,c)$ are listed in Table~\ref{table3}, and the 95\%
bootstrap confidence
intervals of the estimates are also provided. The uninfected cell
proliferation rate ($\lambda$) was
estimated as 46.52 cells per day, the average number of virions
produced by one infected cell ($N$)
was estimated as 1300 per day and the clearance rate of free virions
was 4.35 per day which
corresponds to a half-life of 3.8 hours. All these estimates are in the
ballpark of similar
estimates from other methods [Perelson et al. (\citeyear{r46}, \citeyear{r44})]. In Figure
\ref{figure2}(c), the estimated trajectory of
the time-varying parameter $\eta(t)$ (the viral infection rate), is
plotted with 95\% bootstrap
quantile confidence intervals, which shows an initial fluctuation but
converges to a constant after
2 to 3 months.

\section{Discussion}\label{sec6}

In this paper, we have systematically studied numerical solution-based
NLS estimators for general
nonlinear ODE models which the closed-form solutions are not available.
Both constant and
time-varying parameters are considered. For the model involved
time-varying parameters, we
formulated the estimator under the framework of sieve approach. Our
main contribution is the
establishment of the asymptotic properties for the proposed numerical
solution-based NLS estimators
(including the sieve NLS estimator for the time-varying parameter) with
consideration of both
numerical error and measurement error. Our results show that if the
maximum step size of the
$p$-order numerical algorithm goes to zero at a rate faster than
$n^{-1/(p\wedge4)}$, the numerical
error is negligible compared with the measurement error. This provides
guidance in selecting the
step size for numerical evaluations of ODEs. Moreover, we have shown
that the numerical
solution-based NLS estimator and the sieve NLS estimator for the model
with a time-varying parameter
are strongly consistent. The sieve estimator of constant parameters is
asymptotically normal with
the same asymptotic co-variance as that of the case where the true
solution is exactly known, while
the estimator of the time-varying parameter has an optimal convergence
rate under some regularity
conditions. We also obtained the theoretical results for the case when
the step size of the ODE
numerical solver does not go to zero fast enough or the numerical error
is comparable to the
measurement error [see case (ii) of Theorem \ref{Theorem32} and Remark \ref
{Remark3}]. To our
best knowledge, this is the
first time that the sieve method has been extended to the case of ODE
models which have no
closed-form solutions, and the sieve-based theories were used to
establish the asymptotic results
and construct confidence intervals (bands) for both constant and
time-varying parameters. Note that
we only considered a single time-varying parameter in the model, but
the methodologies can be
extended to multiple time-varying parameters although it is more
tedious to implement.

Note that the NLS estimators have good properties under some
assumptions and are more accurate
compared to other estimates such as those proposed in \citet{r53}, \citet{r10} and
\citet{r35}. But the price that we have to pay is the high
computational cost to obtain the
NLS estimates. To reduce the computational burden, we may use the rough
estimates from other methods
[\citet{r53}, \citet{r10}, \citet{r35}] to narrow
down the search range for the
NLS optimization algorithm. More efficient optimization algorithms may
also be employed to speed up
the computation. We are also considering to parallel our global
optimization algorithms on
high-performance computers. Hopefully these efforts can help us to
handle a reasonable size of ODE
models.

This article only considered the initial value problem (IVP), that is,
the initial conditions are
assumed to be given. In practice, the initial conditions can be
estimated from the data. However,
the generalizations of the theoretical results to the cases of
estimated initial conditions and
other boundary value problems as well as constraints on parameters are
not trivial. Also note that,
if there is more than one time-varying parameter in the model, similar
identifiability techniques
in Section~\ref{sec2} may be applied to these parameters one by one,
sequentially. Spline approximation to
these multiple time-varying parameters can be used for estimation. But
the computation and
theoretical results are more complicated in this case. However, these
generalizations are worth
further investigations in future.

\begin{appendix}\label{app}

\section*{Appendix: Proofs}

\begin{Lemma}\label{Lemma1}
Under conditions \textup{A1--A5}, $\sup_{t\in
I}\|\tilde{\mathbf X}(t,\bolds\beta)-{\mathbf X}(t,\beta)\|_{\infty
}=\break O(h^{p\wedge4})$ for any given
$\bolds\beta\in\mathcal{B}$ in (\ref{Const_ODE}).
\end{Lemma}
\begin{pf}
By Theorem 3.4 in Hairer, N{\o}rsett and Wanner [(\citeyear{r18}),
page 160] under conditions A1--A5, for the $p$th order numerical
algorithm (\ref{num_meth}) for (\ref{Const_ODE}), its global
discretization error satisfies
\[
{\max_{0\leq i\leq m-1}}\|\tilde{\mathbf X}(s_i,\bolds{\beta
})-\mathbf{X}(s_i,\bolds{\beta})\|_{\infty} =O(h^p)\qquad
\mbox{for given $\bolds\beta\in\mathcal{B}$}.
\]
When $t$ is not coincident with the grid points of the numerical
algorithm, the cubic Hermite interpolation [\citet{r13}, page 51]
will be used to obtain the solution at time $t$. In this case,
\[
{\sup_{t\in I\setminus\{s_i\dvtx0\leq i\leq m-1\}}}\|\tilde{\mathbf
X}(t,\bolds{\beta})-\mathbf{X}(t,\bolds{\beta})\|_{\infty}=O(h^4).
\]
Then
it follows that
\begin{eqnarray*}
&&\sup_{t\in I}\|\tilde{\mathbf X}(t,\bolds{\beta})-\mathbf
{X}(t,\bolds{\beta})\|_{\infty}\\
&&\qquad\leq
{\sup_{t\in I\setminus\{s_i\dvtx0\leq i\leq m-1\}}}\|\tilde{\mathbf
X}(t,\bolds{\beta})-\mathbf{X}(t,\bolds{\beta})\|_{\infty}\\
&&\qquad\quad{} +{\max_{t\in\{s_i\dvtx0\leq i\leq m-1\}}}\|\tilde{\mathbf
X}(t,\bolds{\beta})-\mathbf{X}(t,\bolds{\beta})\|_{\infty}\\
&&\qquad=O(h^4)+O(h^p).
\end{eqnarray*}
In general, $h$ is less than 1, $O(h^4)+O(h^p)=O(h^{p\wedge4})$,
which completes the proof.
\end{pf}

Moreover, Lemma \ref{Lemma1} can be extended to the ODE model
(\ref{Time_ODE}) with both constant and time-varying parameters, since
for this model, it can be verified that the result of Theorem 3.1 in
Hairer, N{\o}rsett and Wanner [(\citeyear{r18}), page 157] is still valid
for any given $\bolds\beta\in\mathcal{B}$ and $\eta\in\mathcal{A}$
under condition B2 (it can be derived using the Taylor expansion and
the Chain rule), which leads to the same conclusion as Theorem 3.4 in
Hairer, N{\o}rsett and Wanner [(\citeyear{r18}), page 160]. For Theorems
\ref{Theorem31} and \ref{Theorem32}, the proofs for the univariate and
multivariate cases are the same. For presentation and notation
simplicity, we only outline the proof for the univariate case below.
\begin{pf*}{Proof of Theorem \ref{Theorem31}}
Denote $\tilde{M}_n(\bolds\beta)=\frac{1}{n}\sum_{i=1}^n[Y(t_i)-\tilde
{X}(t_i,\bolds\beta)]^2$,
$M_n(\bolds\beta)=\frac{1}{n}\sum_{i=1}^n[Y(t_i)-X(t_i,\bolds\beta)]^2$
and $M(\bolds\beta)=[Y(t)-X(t,\bolds\beta)]^2$.

First, we claim that $E_0[M(\bolds\beta)]$ reaches its unique minimum
at $\bolds\beta=\bolds{\beta}_0$. In fact,
\begin{eqnarray*}
E_0[M(\bolds\beta)]&=&E_0[Y(t)-X(t,\bolds\beta)]^2\\
&=&E_0[Y(t)-X(t,\bolds{\beta}_0)+X(t,\bolds{\beta}_0)-X(t,\bolds
\beta)]^2\\
&=&E_0[Y(t)-X(t,\bolds{\beta}_0)]^2+E_t[X(t,\bolds{\beta
}_0)-X(t,\bolds\beta
)]^2\\
&=&E_0[\varepsilon(t)]^2+E_t[X(t,\bolds{\beta}_0)-X(t,\bolds\beta
)]^2\\
&\geq&E_0[\varepsilon(t)]^2=E_0[M(\bolds{\beta}_0)],
\end{eqnarray*}
where the third equality holds because the intersection term equals
zero according to the following calculation:
\begin{eqnarray*}
& &E_0[\varepsilon(t)][X(t,\bolds{\beta}_0)-X(t,\bolds\beta)]\\
&&\qquad=E_tE_0 \{[\varepsilon(t)][X(t,\bolds{\beta}_0)-X(t,\bolds\beta
)]|t \}\\
&&\qquad=E_t \{[X(t,\bolds{\beta}_0)-X(t,\bolds\beta)]E_0[\varepsilon
(t)] \}\\
&&\qquad=0,
\end{eqnarray*}
because of $E_0[\varepsilon(t)]=0$. Moreover,
$E_t[X(t,\bolds{\beta})-X(t,\bolds{\beta}_0)]^2=0$ if and only if
$\bolds{\beta}=\bolds{\beta}_0$ from assumption A6. Thus the above claim
holds. Under assumption A10, it follows\vspace*{-1pt} that the first-order
derivative $\frac{\partial E_0 [M(\bolds\beta)]}{\partial\bolds
\beta}$
of $E_0 [M(\bolds\beta)]$ at $\bolds{\beta}_0$ equals to zero and the
second-order derivative $\frac{\partial^2 E_0
[M(\bolds\beta)]}{\partial\bolds\beta\,\partial\bolds\beta^T}$ of $E_0
[M(\bolds\beta)]$ at $\bolds{\beta}_0$ is positive definite. By
assumptions A7 and A9, the second-order derivative of $E_0
[M(\bolds\beta)]$ in a small neighborhood of $\bolds{\beta}_0$ is
bounded away from 0 and $\infty$. Then the second-order Taylor
expansion of $E_0[M(\bolds\beta)]$ gives that there exists a constant
$0<C<\infty$ such that
\[
E_0[M(\hat{\bolds{\beta}}_n)-M(\bolds{\beta}_0)]
\geq C \|\hat{\bolds{\beta}}_n-\bolds{\beta}_0\|^2.
\]
Thus it is sufficient to prove
$E_0[M(\bolds{\beta}_0)]-E_0[M(\hat{\bolds{\beta}}_n)]\rightarrow0$,
a.s.

Let $N_1(\varepsilon,\mathcal{Q},\mathcal{F})$ be the covering number
of the class $\mathcal{F}$ in the
probability measure $\mathcal{Q}$, as given in Pollard (\citeyear{r47}, page 25).
From Lemma 4.1 in \citet{r49},
we have that $N_1(\varepsilon,L_2,\mathcal{B})\leq(\frac{3
R_{\bolds\beta
}}{\varepsilon})^{d}$. Let $\mathcal{F}_n$ be the set $\{M_n(\bolds
\beta
)\dvtx\bolds\beta\in\mathcal{B}\}$. With the Taylor expansion, for any
$\bolds{\beta}_1$, $\bolds{\beta}_2\in\mathcal{B}$, we can easily obtain
\[
|M_n(\bolds{\beta}_1)-M_n(\bolds{\beta}_2)|\leq
C\|\bolds{\beta}_1-\bolds{\beta}_2\|,
\]
where $C$ is some constant. Then for any probability measure $Q$, we
have
\[
\sup_Q N_1(\varepsilon,Q,\mathcal{F}_n)\leq N_1(\varepsilon
/C,L_2,\mathcal{B})\leq C\biggl(\frac{1}{\varepsilon}\biggr)^d\qquad \mbox{for
$0<\varepsilon<1$}.
\]
Then by Theorem II.37 in \citet{r47},
$\sup_{\bolds\beta}|M_n(\beta)-E_0 M(\bolds\beta)|\rightarrow0$, a.s.,
under $P_{\bolds{\beta}_0}$. Then we have
$M_n(\hat{\bolds{\beta}}_n)-E_0 [M(\hat{\beta}_n)]\rightarrow0$
and $
M_n(\bolds{\beta}_0)-E_0 [M(\bolds{\beta}_0)]\rightarrow0$, a.s.

Next, by Lemma \ref{Lemma1},
%
%
\begin{eqnarray}\label{gap1}
\tilde{M}_n(\bolds\beta)&=&\frac{1}{n}\sum_{i=1}^n[Y(t_i)-\tilde
{X}(t_i,\bolds\beta)]^2\nonumber\\
&=&\frac{1}{n}\sum_{i=1}^n\bigl[Y(t_i)-X(t_i,\bolds\beta)+O\bigl(n^{-\lambda
(p\wedge4)}\bigr)\bigr]^2\nonumber\\[-8pt]\\[-8pt]
&=&\frac{1}{n}\sum_{i=1}^n[Y(t_i)-X(t_i,\bolds\beta
)]^2+O\bigl(n^{-\lambda
(p\wedge4)}\bigr)\nonumber\\
&=&M_n(\bolds\beta)+O\bigl(n^{-\lambda(p\wedge4)}\bigr).\nonumber
\end{eqnarray}
Then
\begin{eqnarray*}
& &\tilde{M}_n(\hat{\bolds{\beta}}_n)-E_0 [M(\bolds{\beta}_0)]\\
&&\qquad\geq\tilde{M}_n(\hat{\bolds{\beta}}_n)-E_0 [M(\hat{\bolds
{\beta}}_n)]\\
&&\qquad=M_n (\hat{\bolds{\beta}}_n)+O\bigl(n^{-\lambda(p\wedge4)}\bigr)-E_0
[M(\hat{\bolds{\beta}}_n)]
\end{eqnarray*}
and
\begin{eqnarray*}
& &\tilde{M}_n(\hat{\bolds{\beta}}_n)-E_0 [M(\bolds{\beta}_0)]\\
&&\qquad\leq\tilde{M}_n(\bolds{\beta}_0)-E_0 [M(\bolds{\beta}_0)]\\
&&\qquad=M_n (\bolds{\beta}_0)+O\bigl(n^{-\lambda(p\wedge4)}\bigr)-E_0
M(\bolds{\beta}_0).
\end{eqnarray*}
Hence $\tilde{M}_n(\hat{\bolds{\beta}}_n)-E_0
[M(\bolds{\beta}_0)]\rightarrow0$, a.s. Thus
\begin{eqnarray*}
& &|E_0[M(\hat{\bolds{\beta}}_n)]-E_0[M(\bolds{\beta}_0)]|\\
&&\qquad\leq |\tilde{M}_n(\hat{\bolds{\beta}}_n)-E_0
[M(\hat{\bolds{\beta}}_n)]|+|\tilde{M}_n(\hat{\bolds{\beta}}_n)-E_0
[M(\bolds{\beta}_0)]|\rightarrow0\qquad \mbox{a.s.}
\end{eqnarray*}
Since $\bolds{\beta}_0$ is the unique minimum point for
$E_0[M(\bolds\beta)]$, $\hat{\bolds{\beta}}_n$ is almost surely
consistent with respect to $P_{\bolds{\beta}_0}$.
\end{pf*}
\begin{pf*}{Proof of Theorem \ref{Theorem32}} For the proof of part (i),
it suffices to verify conditions of Theorem 2 in \citet{r48}.
Denote
$\tilde{G}_n(\bolds\beta)=\frac{1}{n}\sum_{i=1}^n[Y(t_i)-\tilde
{X}(t_i,\bolds\beta)]
\frac{\partial\tilde{X}(t_i,\bolds\beta)}{\partial\bolds\beta}$,
$G_n(\bolds\beta)=\frac{1}{n}\sum_{i=1}^n[Y(t_i)-X(t_i,\bolds\beta
)]\,\frac
{\partial
X(t_i,\bolds\beta)}{\partial\bolds\beta}$ and
$G(\bolds\beta)=\break E_0[Y(t)-X(t,\bolds\beta)]\,\frac{\partial
X(t,\bolds\beta)}{\partial\bolds\beta}$. Obviously,
$\tilde{G}_n(\hat{\bolds{\beta}}_n)=0$ and $
G(\bolds{\beta}_0)=E_tE_0(\{[Y(t)-X(t,\bolds{\beta}_0)] \,\frac
{\partial
X(t,\bolds{\beta}_0)}{\partial\bolds{\beta}_0}\}|t)=0$ from
$E_0[Y(t)|t]=X(t,\bolds{\beta}_0)$.

First, we verify the following result:
$\sqrt{n}[\tilde{G}_n(\bolds{\beta}_0)-G(\bolds{\beta
}_0)]\stackrel
{d}{\rightarrow}
N(0,H_1)$. For fixed $t$, according to the multivariate inequality of
Kolmogorov type for $L_2$-norms of derivatives [Babenko, Kofanov and
Pichugov (\citeyear{r3}),
page 9], we have $\|\frac{\partial\tilde{\mathbf X}(t,\bolds\beta
)}{\partial\bolds\beta}-\frac{\partial{\mathbf X}(t,\bolds\beta
)}{\partial\bolds\beta}\| \leq C
\|\frac{\partial^2\tilde{\mathbf X}(t,\bolds\beta)}{\partial
\bolds\beta\,\partial\bolds\beta^T}-\frac{\partial^2
{\mathbf X}(t,\bolds\beta)}{\partial\bolds\beta\,\partial\bolds
\beta^T}\|_{\infty}^{1/2}
\|\tilde{\mathbf X}(t,\bolds\beta)-{\mathbf X}(t,\bolds\beta)\|
_{\infty}^{1/2}\leq C'\|\tilde{\mathbf X}(t,\bolds\beta)-{\mathbf
X}(t,\bolds\beta)\|_{\infty}^{1/2}$ for two
constants $C$ and $C'$, where the second inequality holds because of
the uniform boundedness of both $\frac{\partial^2 {\mathbf
X}(t,\bolds\beta)}{\partial\bolds\beta\,\partial\bolds\beta^T}$ and
$\frac{\partial^2\tilde{\mathbf X}(t,\bolds\beta)}{\partial\bolds
\beta\,\partial\bolds\beta^T}$ under
conditions A7 and A8. Based on $\sup_{t\in I}\|\tilde{\mathbf
X}(t,\bolds\beta)-{\mathbf X}(t,\beta)\|_{\infty}=O(n^{-\lambda
(p\wedge
4)})$ from Lemma \ref{Lemma1}, it follows that $\|\frac{\partial
\tilde{\mathbf X}(t,\bolds\beta)}{\partial\bolds\beta}-\frac
{\partial{\mathbf X}(t,\bolds\beta)}{\partial\bolds\beta}\|
=O(n^{-\lambda(p\wedge4)/2})$.
Considering that $Y(t_i)-X(t_i,\bolds{\beta}_0)$ and
$\frac{\partial\tilde{X}(t_i,\bolds{\beta}_0)}{\partial\bolds
\beta_0}$
are bounded, we have
\begin{eqnarray*}
\hspace*{-5pt}& &\sqrt{n}[\tilde{G}_n(\bolds{\beta}_0)-G(\bolds{\beta}_0)]\\
\hspace*{-5pt}&&\qquad=\frac{1}{\sqrt{n}}\sum_{i=1}^n[Y(t_i)-\tilde{X}(t_i,\bolds
{\beta
}_0)]\,\frac{\partial\tilde{X}(t_i,\bolds{\beta}_0)}{\partial\bolds
\beta_0}\\
\hspace*{-5pt}&&\qquad=\frac{1}{\sqrt{n}}\sum_{i=1}^n\bigl[Y(t_i)-X(t_i,\bolds{\beta
}_0)+O\bigl(n^{-\lambda(p\wedge
4)}\bigr)\bigr]\biggl[\frac{\partial X(t_i,\bolds{\beta}_0)}{\partial\bolds{\beta}_0}
+O\bigl(n^{-\lambda(p\wedge4)/2}\bigr)\biggr]\\
\hspace*{-5pt}&&\qquad=\frac{1}{\sqrt{n}}\sum_{i=1}^n[Y(t_i)-X(t_i,\bolds{\beta
}_0)]\,\frac
{\partial
X(t_i,\bolds{\beta}_0)}{\partial\bolds{\beta}_0}+O\bigl(n^{-\lambda
(p\wedge
4)/2+1/2}\bigr).
\end{eqnarray*}
When $\lambda>1/(p\wedge4)$, $O(n^{-\lambda(p\wedge
4)/2+1/2})=o(1)$. So for the above expression, we have
\begin{eqnarray*}
& &\sqrt{n}[\tilde{G}_n(\bolds{\beta}_0)-G(\bolds{\beta}_0)]\\
&&\qquad=\frac{1}{\sqrt{n}}\sum_{i=1}^n[Y(t_i)-X(t_i,\bolds{\beta
}_0)]\,\frac
{\partial
X(t_i,\bolds{\beta}_0)}{\partial\bolds{\beta}_0}+o(1)\\
&&\qquad=\sqrt{n}[G_n(\bolds{\beta}_0)-G(\bolds{\beta}_0)]+o(1).
\end{eqnarray*}
Based on the general central limit theorem,
$\sqrt{n}[G_n(\bolds{\beta}_0)-G(\bolds{\beta}_0)]\rightarrow
N(0,\mathbf{H}_1)$ with
\[
\mathbf{H}_1=E_0[Y(t)-X(t,\bolds{\beta}_0)]^2\biggl[\frac{\partial
X(t,\bolds{\beta}_0)}{\partial\bolds{\beta}_0}\biggr]^{\otimes2}=\sigma_0^2
E_t\biggl[\frac{\partial
X(t,\bolds{\beta}_0)}{\partial\bolds{\beta}_0}\biggr]^{\otimes2}.
\]

Second, let $\delta_n\downarrow0$. For
$\|\bolds\beta-\bolds{\beta}_0\|\leq\delta_n$, we want to show that
\[
\sqrt{n}[\tilde{G}_n(\bolds\beta)-G(\bolds\beta)]-\sqrt
{n}[\tilde{G}_n(\bolds
{\beta}_0)-G(\bolds{\beta}_0)]=o_p(1).
\]
In fact, from the first step above, for any $\bolds\beta\in\mathcal{B}$,
we have that
$\sqrt{n}[\tilde{G}_n(\bolds\beta)-G_n(\bolds\beta)]=o_p(1)$. Then
\begin{eqnarray*}
&&\sqrt{n}[\tilde{G}_n(\bolds\beta)-G(\bolds\beta)]-\sqrt
{n}[\tilde
{G}_n(\bolds{\beta}_0)-G(\bolds{\beta}_0)]\\
&&\qquad=\sqrt{n}[G_n(\bolds\beta)-G(\bolds\beta)]-\sqrt{n}[G_n(\bolds
{\beta
}_0)-G(\bolds{\beta}_0)]+o_p(1).
\end{eqnarray*}
From Lemma 4.1 in \citet{r49}, we have that
$N_1(\varepsilon,L_2,\mathcal{B})\leq(\frac{3 R}{\varepsilon
})^{d}$. Let
$\Lambda_n$ be the set $\{G_n(\bolds\beta)\dvtx\bolds\beta\in
\mathcal{B}\}$
for any $X\in\mathcal{X}$. Using a Taylor series expansion, for any
$\bolds{\beta}_1$, $\bolds{\beta}_2\in\mathcal{B}$, we can easily obtain
\[
|G_n(\bolds{\beta}_1)-G_n(\bolds{\beta}_2)|\leq
C\|\bolds{\beta}_1-\bolds{\beta}_2\|,
\]
where $C$ is some constant. Then for any probability measure $Q$, we
have
\[
N_1(\varepsilon,L_2(Q),\Lambda_n)\leq N_1(\varepsilon/C,L_2,\mathcal
{B})\leq
C\biggl(\frac{1}{\varepsilon}\biggr)^d,
\]
and thus
\[
\log N_1(\varepsilon,L_2(Q),\Lambda_n)\leq d\log{\frac
{1}{\varepsilon}}.
\]
Since $\int_0^1\log(1/\varepsilon)\,d\varepsilon<\infty$, $\Lambda_n$
is a
P-Donsker class by Theorem 2.5.2 in \citet{r66}.
Hence
$\sqrt{n}[G_n(\bolds\beta)-G(\bolds\beta)]-\sqrt{n}[G_n(\bolds
{\beta
}_0)-G(\bolds{\beta}_0)]=o_p(1)$.

Third, with some simple calculations, we have
$G(\bolds\beta)=E_t[X(t,\bolds{\beta}_0)-X(t$,\break$\bolds\beta)]\,\frac
{\partial
X(t,\bolds\beta)}{\partial\bolds\beta}$, then
\[
\frac{\partial G(\bolds\beta)}{\partial\bolds\beta}=-\int
\biggl\{\frac{\partial X(t,\bolds\beta)}{\partial\bolds\beta}\biggr\}
^{\otimes2}\,
d\Phi(t)+ \int[X(t,\bolds{\beta}_0)-X(t,\bolds\beta)]\,\frac
{\partial^2
X(t,\bolds\beta)}{\partial\bolds\beta\,\partial\bolds\beta^T}\,
d\Phi(t)
\]
and $\frac{\partial
G(\bolds\beta)}{\partial\bolds\beta}|_{\bolds\beta=\bolds{\beta
}_0}=-E_t
\{\frac{\partial
X(t,\bolds{\beta}_0)}{\partial\bolds{\beta}_0}\}^{\otimes2}$. Denote
$\mathbf{H}_2=E_t \{\frac{\partial
X(t,\bolds{\beta}_0)}{\partial\bolds{\beta}_0}\}^{\otimes2}$.
Then by
using the Taylor series expansion again, the function $G(\bolds\beta)$
is\break Fr\'{e}chet-differentiable at $\bolds\beta_0$ with
nonsingular derivative $\mathbf{H}_2$.

Thus all conditions of Theorem 2 in \citet{r48} are satisfied, then
Theorem~\ref{Theorem32}(i) holds with
$\mathbf{V}_1=\mathbf{H}_2^{-1}\mathbf{H}_1(\mathbf
{H}_2^{-1})^T=\sigma_0^2\{
E_t[\frac{\partial
X(t,\bolds{\beta}_0)}{\partial\bolds{\beta}_0}]^{\otimes2}\}^{-1}$.

For the proof of case (ii) of Theorem \ref{Theorem32}, it is easy to verify the
conditions of Theorem 2 in
\citet{r48} for the asymptotic normality. Now we just need to show
$\tilde{\bolds\beta}=\bolds\beta_0+O(h^{(p\wedge4)/2})$ and
$\tilde
{V}_1=V_1+O(h^{(p\wedge4)/2})$.
Denote $\tilde{M}(\bolds\beta)=[Y(t)-\tilde{X}(t,\bolds\beta)]^2$ and
$\tilde{G}(\bolds\beta)=E_0[Y(t)-\tilde{X}(t,\beta)]\,\frac
{\partial\tilde
{X}(t,\bolds\beta)}{\partial\bolds\beta}$.
Since $E_0[\tilde{M}(\bolds\beta)]$ reaches its minimum at $\bolds
\beta=
\tilde{\bolds\beta}$, then the
first-order derivative of $E_0[\tilde{M}(\bolds\beta)]$ at $\tilde{
\bolds\beta}$ equals 0, that is,
$\tilde{G}(\tilde{\bolds\beta})=0$. Then similar to the proof of
case (i)
above, we have
\begin{eqnarray*}
\tilde{G}(\bolds\beta)&=&E_0[Y(t)-\tilde{X}(t,\bolds\beta)]\,\frac
{\partial
\tilde{X}(t,\bolds\beta)}{\partial\bolds\beta}\\
&=& E_0[Y(t)-X(t,\bolds\beta)+O(h^{p\wedge4})]\biggl[\frac{\partial
X(t,\bolds
\beta)}{\partial\bolds\beta}+O\bigl(h^{(p\wedge4)/2}\bigr)\biggr]\\
&=&G(\bolds\beta)+O\bigl(h^{(p\wedge4)/2}\bigr).
\end{eqnarray*}
It follows that
$\tilde{G}(\tilde{\bolds\beta})=G(\tilde{\bolds\beta
})+O(h^{(p\wedge
4)/2})$, then $G(\tilde{\bolds\beta})=O(h^{(p\wedge4)/2})$ from
$\tilde{G}(\tilde{\bolds\beta})=0$. The Taylor series expansion yields
that there exist constants $0<c_1,c_2<\infty$ such that
\[
c_1\|\tilde{\bolds\beta}-\bolds\beta_0\|\leq
|G(\tilde{\bolds\beta})-G(\bolds\beta_0)|\leq
c_2\|\tilde{\bolds\beta}-\bolds\beta_0\|.
\]
Thus $\|\tilde{\bolds\beta}-\bolds\beta_0\|=O(h^{(p\wedge4)/2})$ from
$G(\bolds\beta_0)=0$. Similarly we can show that
$\|\tilde{\mathbf{V}}_1-\mathbf{V}_1\|=O(h^{(p\wedge4)/2})$.
\end{pf*}

Some definitions and notation are necessary in order to prove
Theorems \ref{Theorem41}--\ref{Theorem43}. Denote
$\tilde{M}_n(\bolds\theta)=\frac1n\sum_{i=1}^n\sum_{j=1}^K[Y_j(t_i)-
\tilde{\mathbf{X}}_j(t_i,\bolds\beta,\eta(t_i))]^2$,
$M_n(\bolds\theta)=\frac1n\sum_{i=1}^n\sum_{j=1}^K[Y_j(t_i)-
\mathbf{X}_j(t_i,\bolds\beta,\eta(t_i))]^2$
and
$M(\bolds\theta)=\sum_{j=1}^K[Y_j-\mathbf{X}_j(t,\bolds\beta,\eta
(t))]^2$.
We define a semidistance $\rho$ on $\Theta$ as
\[
\rho^2(\bolds\theta,\bolds{\theta}_0)=E_0\{(\bolds\beta-\bolds
{\beta}_0)^T\dot
{M}_1(\bolds\theta)+\dot{M}_2(\bolds\theta)[\eta-\eta_0]\}^2,
\]
where $\dot{M}_1$ is the score function of $M$ for $\bolds\beta$, and
$\dot{M}_2$ is the score operator of $M$ for~$\eta$, both evaluated
at the true parameter value $\bolds{\theta}_0$. Similarly to the proof
in Huang and Rossini [(\citeyear{r25}), page 966] when $\mathbf{V}_2(\bolds
{\theta}_0)$,
defined in assumption B6, is positive definite, and $\dot{M}_1$ and
$\dot{M}_2$ are bounded away from $+\infty$ and $-\infty$, if
$\rho(\hat{\bolds{\theta}}_n,\bolds{\theta}_0) =O_p (r_n)$, then
$d(\hat{\bolds{\theta}}_n,\bolds{\theta}_0)=O_p (r_n)$; and if
$\rho(\hat{\bolds{\theta}}_n,\bolds{\theta}_0)\rightarrow0$ almost
surely under $P_{\bolds{\theta}_0}$, then
$d(\hat{\bolds{\theta}}_n,\bolds{\theta}_0)\rightarrow0$ almost surely
under $P_{\bolds{\theta}_0}$.
\begin{pf*}{Proof of Theorem \ref{Theorem41}}
Similarly to the proof of
Theorem \ref{Theorem31}, we have that $E_0[M(\bolds{\theta}_0)]$ reaches its
unique minimum at $\bolds\theta=\bolds\theta_0$. It follows that
\[
E_0[M(\bolds{\theta}_0)-M(\hat{\bolds{\theta}}_n)]
\geq C \rho^2(\hat{\bolds{\theta}}_n,\bolds{\theta}_0),
\]
where $C$ is some constant. Thus if
$E_0[M(\bolds{\theta}_0)]-E_0[M(\hat{\bolds{\theta
}}_n)]\rightarrow0$,
almost surely under $P_{\bolds{\theta}_0}$, then
$d(\hat{\bolds{\theta}}_n,\bolds{\theta}_0)\rightarrow0$, almost surely
under $P_{\bolds{\theta}_0}$.\vspace*{1pt}

Let $\mathcal{A}_n^{\delta}$ be the set $\{\eta\in\mathcal{A}_n,
\|\eta-\eta_{n0}\|_2\leq\delta\}$ and
$N_2(\varepsilon,L_{\infty},\mathcal{A}_n^{\delta})$ be its bracketing
number with respect to $L_{\infty}$ [see Definition 2.1.6, \citet{r66}], where $\eta_{n0}$ is the map point of
$\eta_0$ in the sieve $\mathcal{A}_n$.
By the calculation of Shen and Wong [(\citeyear{r59}), page 597] for any
$\varepsilon\leq\delta$,
we have
\[
N_2(\varepsilon,L_{\infty},\mathcal{A}_n^{\delta})\leq C(\delta
/\varepsilon)^{N},
\]
where $N=q+l$ is the number of B-splines basis functions. Let $\mathcal
{F}_n$ be the set $\{M_n(\bolds\theta)\dvtx\|\bolds\beta-\bolds
{\beta}_0\|\leq
\delta,
\eta\in\mathcal{A}_n, \|\eta-\eta_{n0}\|_2\leq\delta\}$.
For any $\bolds{\theta}_1$, $\bolds{\theta}_2\in\Theta_n$, we can easily
obtain
\[
|M_n(\bolds{\theta}_1)-M_n(\bolds{\theta}_2)|\leq
C(\|\bolds{\beta}_1-\bolds{\beta}_2\|+\|\eta_1-\eta_2\|_{\infty})
\]
using Taylor's expansion. Hence
\begin{eqnarray*}
N_2(\varepsilon,L_{\infty},\mathcal{F}_n)&\leq& N_1(\varepsilon
/2,L_2,\mathcal{B})\times N_2(\varepsilon/2,L_{\infty},\mathcal
{A}_n^{\delta})\\
&\leq& C(3R_d/\varepsilon)^d(\delta/\varepsilon)^{N}\\
&\leq& C'(1/\varepsilon)^{N+d}.
\end{eqnarray*}
Note that, since $N_2(\varepsilon,L_{\infty},\mathcal{F}_n)$ depends on
$n$ in the above expression, we cannot directly use Theorem II.37 in
\citet{r47} to obtain ${\sup_{\mathcal{F}_n}}|M_n(\bolds\theta)-E_0
[M(\bolds\theta)]|\rightarrow0$, a.s., under $P_{\bolds{\theta}_0}$.
Fortunately, we can still get this result based on (A.2) in
\citet{r76}. Thus we have $M_n(\hat{\bolds{\theta}}_n)-E_0
[M(\hat{\bolds{\theta}}_n)]\rightarrow0$ and $
M_n(\bolds{\theta}_{n0})-E_0 [M(\bolds{\theta}_{n0})]\rightarrow0$,
a.s., where $\bolds{\theta}_{n0}$ is the map point of $\bolds{\theta}_0$
in the sieve~$\Theta_n$.

From the extension of Lemma \ref{Lemma1} for any given $\bolds\beta
\in\mathcal{B}$
and $\eta(t)\in\mathcal{A}$ in (\ref{Time_ODE}), similarly to
(\ref{gap1}), we have
%
%
\begin{equation}\label{gap2}
\tilde{M}_n(\bolds\theta) =M_n(\bolds\theta)+O\bigl(n^{-\lambda
(p\wedge4)}\bigr).
\end{equation}
Then the remaining steps are similar to those in the proof of Theorem
\ref{Theorem31}.
\end{pf*}
\begin{pf*}{Proof of Theorem \ref{Theorem42}} We apply Theorem 3.4.1 in
\citet{r66} to obtain the rate of convergence.

For $\bolds\theta_{n0}$ in the proof of Theorem \ref{Theorem41}, define
$\bolds\theta_{n0}\mapsto\rho_1(\bolds\theta,\bolds\theta _{n0})$ be a
map from $\Theta_n$ to $[0,\infty)$ as
$\rho_1^2(\bolds\theta,\bolds\theta_{n0})=E_0 [M(\bolds\theta)]-E_0
[M(\bolds\theta_{n0})]$. Choose
$\delta_n=\rho(\bolds\theta_0,\bolds\theta_{n0})$. For
$\delta_n<\delta<\infty$, denote
$\Omega=\{\bolds\theta\dvtx\bolds\theta\in\Theta_n,
\delta/2<\rho(\bolds\theta,\bolds\theta_{n0})\leq\delta\}$. From the
definition of~$\rho_1$, we have $\sup_{\Omega}E_0
[M(\bolds\theta_{n0})]-E_0 [M(\bolds\theta)]\leq-\frac{\delta^2}{4}$.

Let $\Xi_n$ be the set
$\{M_n(\bolds\theta)-M(\bolds\theta_{n0})\dvtx\bolds\theta\in
\Theta_n\}$ and $\tilde{J}(\delta,L_2(P),\Xi_n)$ be the $L_2(P)$-norm
bracketing integral of the sieve $\Theta_n$. From the proof of Theorem
\ref{Theorem41}, we have
\begin{eqnarray*}
\tilde{J}(\delta,L_2(P),\Xi_n)&=&\int_0^{\delta}\sqrt{1+\log
N_2(\varepsilon,L_2(P),\Xi_n)}
\,d\varepsilon\\
&\leq& \int_0^{\delta}\sqrt{1+\log N_2(\varepsilon,L_{\infty},\Xi_n)}
\,d\varepsilon\\
&\leq& C N^{1/2}\delta.
\end{eqnarray*}
Let
\[
\phi_n(\delta)=\tilde{J}(\delta,L_2(P),\Xi_n)\biggl(1+\frac{\tilde
{J}(\delta,
L_2(P),\Xi_n)}{\delta^2\sqrt{n}}\biggr)=N^{1/2}\delta+\frac
{N}{\sqrt{n}}.
\]
Obviously, $\phi_n(\delta)/\delta^{1+\tau}$ is a decreasing
function in
$\delta$ for $0<\tau<1$.
Then by Lemma 3.4.2 in \citet{r66}, we have
\[
E_0\Bigl[\sup_{\Omega}\sqrt{n}(M_n-M)(\bolds\theta-\bolds\theta
_{n0})\Bigr]\preceq
\phi_n(\delta).
\]
For $\lambda>1/[2(p\wedge4)]$, from (\ref{gap2}), it follows that
\[
\sqrt{n}[\tilde{M}_n(\bolds\theta)-M_n(\bolds\theta)]=O\bigl(n^{
1/2-\lambda
(p\wedge
4)}\bigr)=o(1).
\]
Then we have that
\[
E_0\Bigl[\sup_{\Omega}\sqrt{n}(\tilde{M}_n-M)(\bolds\theta-\bolds
\theta
_{n0})\Bigr]\preceq
\phi_n(\delta).
\]
Then the conditions of Theorem 3.4.1 in \citet{r66}
are satisfied for the $\delta_n$, $\rho_1$ and
$\phi_n(\delta)$ above. Therefore we have
$-r_n^2\rho_1(\hat{\bolds{\theta}}_n,\bolds\theta_{n0})=O_p(1)$, where
$r_n$ satisfies $r_n^2 \phi_n(\frac{1}{r_n})\leq\sqrt{n}$. It
follows that $r_n=N^{-1/2}n^{1/2}=n^{(1-v)/2}$. Thus
$\rho_1(\hat{\bolds{\theta}}_n,\bolds\theta_{n0})=O_p(n^{-(1-v)/2})$.

Now, we define a distance $\rho_2$ as
\[
\rho_2(\bolds\theta_1,\bolds\theta_2)=\|\bolds\beta_1-\bolds
\beta_2\|+\|\eta
_1-\eta_2\|_{\infty}.
\]
Let $\varsigma$ be a positive constant. Similarly to the proof of
Theorem 3.2 in \citet{r24}, it is easy to follow that for any
$\bolds\theta$ with $\rho_2(\bolds\theta,\bolds\theta_{n0})\leq
\varsigma$,
there exist constants $0<c_1,c_2<\infty$ such that
\[
-c_1d^2(\bolds\theta,\bolds\theta_{n0})+O_p(n^{-2v\varrho})\leq
-\rho^2_1(\bolds\theta,\bolds\theta_{n0})\leq
-c_2d^2(\bolds\theta,\bolds\theta_{n0})+O_p(n^{-2v\varrho}).
\]
Therefore, for a constant $c_2>0$,
\[
c_2d^2(\hat{\bolds\theta}_n,\bolds\theta_{n0})\leq
O_p\bigl(n^{-2v\varrho}+n^{-(1-v)}\bigr).
\]
Because $d(\bolds\theta_{n0},\bolds\theta_0)\leq
\rho_2(\bolds\theta_0,\bolds\theta_0)=O_p(n^{-v\varrho})$, we have
$d(\hat{\bolds\theta}_n,\bolds\theta_0)=O_p(n^{-v\varrho}+n^{-(1-v)/2})$.
\end{pf*}
\begin{pf*}{Proof of Theorem \ref{Theorem43}}
We prove this theorem using Theorem 6.1 in \citet{r71}. It
suffices to validate conditions A1--A6 of Theorem~6.1 in
\citet{r71}. From the proof of Theorems \ref{Theorem41} and
\ref{Theorem42} above, it is easy to see that condition A1 regarding
consistency and rate of convergence and condition A2 for Theorem 6.1 in
\citet{r71} hold.

For condition A3, we need to calculate the pseudo-information
matrix. For any fixed $\eta\in\mathcal{A}$, let $\mathcal{A}_0=\{
\eta
_{\omega}(\cdot)\dvtx\omega$ in a neighborhood of $0\in
\mathcal{R}\}$ be a smooth curve in $\mathcal{A}$ running through
$\eta_0$
at $\omega=0$, that is, $\eta_{\omega=0}(t)=\eta_0(t)$. Denote
$\frac{\partial}{\partial\omega}\,\eta_{\omega}(t)|_{\omega=0}=a(t)$
and the space generated by such $a(t)$ as $\Upsilon$. The score
functions of $\bolds\beta$ and $\eta$ are
\begin{eqnarray*}
\dot{M}_1&=&\frac{\partial
M}{\partial\bolds\beta}=-2\sum_{j=1}^K(Y_j-\mathbf{X}_j)\,\frac
{\partial
\mathbf{X}_j}{\partial\bolds\beta_0},\\
\dot{M}_2[a]&=&\frac{\partial M}{\partial\eta_0}=
-2\sum_{j=1}^K(Y_j-\mathbf{X}_j)\,\frac{\partial
\mathbf{X}_j}{\partial\eta_0}a(t).
\end{eqnarray*}
We also set
\begin{eqnarray*}
\dot{M}_{11}&=&\frac{\partial^2
M}{\partial\bolds\beta_0\,\partial\bolds\beta_0^T}
=2\sum_{j=1}^K\biggl[\frac{\partial\mathbf{X}_j}{\partial\bolds\beta
_0}\,\frac
{\partial\mathbf{X}_j}{\partial
\bolds\beta_0^T}-(Y_j-\mathbf{X}_j)\,\frac{\partial^2
\mathbf{X}_j}{\partial\bolds\beta_0\,\partial\bolds\beta_0^T}\biggr],
\\
\dot{M}_{12}[a]&=&\dot{M}_{21}^T[a]=\frac{\partial^2
M}{\partial\bolds\beta_0\,\partial\bolds\eta_0}=2\sum
_{j=1}^K\biggl[\frac{\partial
\mathbf{X}_j}{\partial\bolds\beta_0}\,\frac{\partial
\mathbf{X}_j}{\partial\eta_0}-(Y_j-\mathbf{X}_j)\,\frac{\partial^2
\mathbf{X}_j}{\partial\bolds\beta_0\,\partial\eta_0}\biggr]a(t)
\end{eqnarray*}
and
\[
\dot{M}_{22}[a_1,a_2]=\frac{\partial^2 M}{\partial
\eta_0^2}=2\sum_{j=1}^K\biggl[\biggl(\frac{\partial
\mathbf{X}_j}{\partial\eta_0}\biggr)^2-(Y_j-\mathbf{X}_j)\,\frac{\partial^2
\mathbf{X}_j}{\partial\eta_0^2}\biggr]a_1(t)a_2(t),
\]
where $a_1(t), a_2(t)\in\Upsilon$. Following the idea from the
proofs of the asymptotic results for semiparametric M-estimator in
\citet{r36} and \citet{r71}, we assume that the
special perturbation direction
$\mathbf{a}^*(t)=(a_1^*(t),\ldots,a_d^*(t))^T$ with
$a_i^*(t)\in\Upsilon$ for $1\leq i\leq d$, satisfies
$E_0\{\dot{M}_{12}[a]-\dot{M}_{22}[\mathbf{a}^*,a]\}=0$ for any
$a\in\Upsilon$. Some calculations yield
\begin{eqnarray*}
&&E_0\{\dot{M}_{12}[a]-\dot{M}_{22}[\mathbf{a}^*,a]\}\\
&&\qquad=2\sum_{j=1}^K E_0\biggl[\biggl\{\frac{\partial
\mathbf{X}_j}{\partial\bolds\beta_0}\,\frac{\partial
\mathbf{X}_j}{\partial\eta_0}-(Y_j-\mathbf{X}_j)\,\frac{\partial^2
\mathbf{X}_j}{\partial\bolds\beta_*\,\partial\eta_0}\biggr\}a(t)\\
&&\qquad\quad\hspace*{39.7pt}{}-\biggl\{\biggl(\frac{\partial
\mathbf{X}_j}{\partial\eta_0}\biggr)^2-(Y_j-\mathbf{X}_j)\,\frac{\partial^2
\mathbf{X}_j}{\partial\eta_0^2}\biggr\}a(t)\mathbf{a}^*(t)\biggr]\\
&&\qquad=2\sum_{j=1}^K E_tE_0 \biggl(\biggl[\biggl\{\frac{\partial
\mathbf{X}_j}{\partial\bolds\beta_0}\frac{\partial
\mathbf{X}_j}{\partial\eta_0}-(Y_j-\mathbf{X}_j)\,\frac{\partial^2
\mathbf{X}_j}{\partial\bolds\beta_*\,\partial\eta_0}\biggr\}a(t)\\
&&\qquad\quad\hspace*{58.6pt}{}-\biggl\{\biggl(\frac{\partial
\mathbf{X}_j}{\partial\eta_0}\biggr)^2-(Y_j-\mathbf{X}_j)\,\frac{\partial^2
\mathbf{X}_j}{\partial\eta_0^2}\biggr\}a(t)\mathbf{a}^*(t)\biggr]\Big|t \biggr).
\end{eqnarray*}
%
It follows that
%
%
\begin{equation}\label{general_at}\qquad
\mathbf{a}^*(t)=\frac{\sum_{j=1}^KE_*[\{{\partial
\mathbf{X}_j}/{\partial\bolds\beta_0}\,{\partial
\mathbf{X}_j}/{\partial\eta_0}-(Y_j-\mathbf{X}_j){\partial^2
\mathbf{X}_j}/{\partial\bolds\beta_0\,\partial\eta_0}\}|t]}{\sum
_{j=1}^KE_0[\{
({\partial
\mathbf{X}_j}/{\partial\eta_0})^2-(Y_j-\mathbf{X}_j){\partial^2
\mathbf{X}_j}/{\partial\eta_0^2}\}|t]}.
\end{equation}
Since
$E_0[\mathbf{Y}(t)|t]=\mathbf{X}(t)$, $\mathbf{a}^*(t)$ in (\ref
{general_at})
can be simplified as
%
%
\begin{equation}\label{simple_at}
\mathbf{a}^*(t)=\frac{\sum_{j=1}^K{\partial
\mathbf{X}_j}/{\partial\bolds\beta_0}\,{\partial
\mathbf{X}_j}/{\partial\eta_0}}{\sum_{j=1}^K({\partial
\mathbf{X}_j}/{\partial\eta_0})^2}.
\end{equation}
%
For $K\geq2$, both
%
%
\begin{equation}\label{S1}\qquad
\mathbf{S}_1=E_0(\dot{M}_{11}-\dot{M}_{12}[\mathbf{a}^*])=2\sum
_{j=1}^KE_t\biggl[\frac{\partial
\mathbf{X}_j}{\partial\bolds\beta_0}\,\frac{\partial
\mathbf{X}_j}{\partial\bolds\beta_0^T}-\frac{\partial
\mathbf{X}_j}{\partial\bolds\beta_0}\,\frac{\partial
\mathbf{X}_j}{\partial\eta_0}\mathbf{a}^*(t)\biggr]
\end{equation}
and
%
%
\begin{equation}\label{S2}
\mathbf{S}_2=E_0(\dot{M}_1-\dot{M}_2[\mathbf{a}^*])^{\otimes
2}=4\sum_{j=1}^K\sigma_j^2E_t\biggl\{\frac{\partial
\mathbf{X}_j}{\partial\bolds\beta_0}-\frac{\partial
\mathbf{X}_j}{\partial\eta_0}\mathbf{a}^*(t)\biggr\}^{\otimes2}
\end{equation}
are nonsingular. Let
$\mathbf{V}_2=\mathbf{S}_1^{-1}\mathbf{S}_2(\mathbf{S}_1^{-1})^T$. Thus
condition A3 of finite variance for Theorem 6.1 in \citet{r71} is
satisfied.


Conditions A4 and A5 for Theorem 6.1 in \citet{r71} can be
verified by similar arguments as condition (i) and C3 in the proof of
Theorem 4 in \citet{r76}, respectively. Condition A6 of
smoothness of the model can be easily verified using a straightforward
Taylor expansion where $n^{-c_1}$ is just the rate of convergence in
Theorem \ref{Theorem42} and faster than $n^{-1/4}$, and $c_2=2$, which
completes the proof.
\end{pf*}
\begin{pf*}{Proof of Proposition \ref{Proposition1}}
Let $\mathcal {G}$ be the set of a real valued functions $g$ on $[a,b]$
which are absolutely continuous and satisfy $\int_a^b g^2(t)\,dt<\infty$
and \mbox{$E_tg(t)=0$}. Similarly to the proof of Theorem \ref{Theorem43}, for
any fixed $\eta\in\mathcal{A}$, let
$\mathcal{A}_0=\{\eta_\omega(\cdot)\dvtx\omega$ in a neighborhood
of $0\in\mathcal{R}\}$ be a smooth curve in $\mathcal{A}$ running
through $\eta_0$ at $\omega=0$, that is,
$\eta_{\omega=0}(t)=\eta_0(t)$. Denote $\frac{\partial}{\partial
\omega}\eta_\omega(t)|_{\omega=0}=a(t)$ and restrict
$a(t)\in\mathcal{G}$. Denote the space generated by such $a(t)$ as
$\Upsilon$. The score functions of $\beta$ and $\eta$ are
\begin{eqnarray*}
\dot{M}_1&=&\frac{\partial M}{\partial\beta}=-2(Y-X)\,\frac{\partial
X}{\partial\xi},\\
\dot{M}_2[a]&=&\frac{\partial M}{\partial\eta}=
-2(Y-X)\,\frac{\partial X}{\partial\xi}a(t)
\end{eqnarray*}
with $\xi=\beta_0+\eta_0(t)$. Let $\dot{P}$ be the linear span of
$\dot{M}_2[a]$. Since $E_0\{\dot{M}_1\dot{M}_2[a]\}=0$ for any
$a(t)\in\Upsilon$, it follows that $\dot{M}_1$ is orthogonal to
$\dot{P}$. Thus the efficient score function of $\beta$ is just
$\dot{M}_1$. Then the pseudo-information is $E_0[\dot{M}_1^2]$. The
rest of the proof is similar to that of Theorem \ref{Theorem43}, where the
efficient score function and the pseudo-information are updated as
discussed before, and the least favorable direction can be selected
by any $a\in\Upsilon$.
\end{pf*}
\end{appendix}

\section*{Acknowledgments}
The authors thank Drs. Hua Liang, Xing Qiu and Jianhua Huang for
helpful discussions, and Ms. Jeanne
Holden-Wiltse for assistance in editing the manuscript. We also highly
appreciate the two referees
and the Associate Editors for their insightful comments and useful suggestions that
have helped us to greatly
improve this manuscript.

\printaddresses


\begin{thebibliography}{99}

\bibitem[\protect\citeauthoryear{Adams}{2005}]{r1}
\textsc{Adams, B. M.} (2005). Non-parametric parameter
estimation and clinical data fitting with a model of HIV infection.
Ph.D. thesis, North Carolina State Univ.
\MR{2623319}

\bibitem[\protect\citeauthoryear{Anderson and May}{1991}]{r2}
\textsc{Anderson, R. M.} and \textsc{May, R. M.} (1991). \textit{Infectious Diseases of
Humans: Dynamics and Control}. Oxford Univ. Press, Oxford.

\bibitem[\protect\citeauthoryear{Audoly et al.}{2001}]{Audoly2001}
\textsc{Audoly, S.}, \textsc{Bellu, G.}, \textsc{D'Angio, L.},
\textsc{Saccomani, M. P.} and \textsc{Cobelli, C.} (2001). Global identifiability
of nonlinear models of biological systems.
\textit{IEEE Trans. Biomed. Eng.} \textbf{48} 55--65.

\bibitem[\protect\citeauthoryear{Babenko, Kofanov and
Pichugov}{1996}]{r3}
\textsc{Babenko, V. F.}, \textsc{Kofanov, V. A.} and
\textsc{Pichugov, S. A.} (1996). Multivariate inequalities of
Kolmogorov type and their applications. In \textit{Multivariate
Approximation and Splines} (G. Nuraberger, J. W. Schmidt and G.
Walz, eds.) 1--12. Birkhauser, Basel.
\MR{1484990}

\bibitem[\protect\citeauthoryear{Bard}{1974}]{r4}
\textsc{Bard, Y.} (1974). \textit{Nonlinear Parameter Estimation}.
Academic Press, New York.
\MR{0326870}

\bibitem[\protect\citeauthoryear{Bellman and {\AA}str\"{o}m}{1970}]{Bellman1970}
\textsc{Bellman, R.} and \textsc{{\AA}str\"{o}m, K. J.} (1970). On
structural identifiability. \textit{Math. Biosci.} \textbf{7} 329--339.

\bibitem[\protect\citeauthoryear{Benson}{1979}]{r5}
\textsc{Benson, M.} (1979). Parameter fitting in dynamic models.
\textit{Ecol. Mod.} \textbf{6} 97--115.

\bibitem[\protect\citeauthoryear{Bickel et al.}{1993}]{r6}
\textsc{Bickel, P. J.}, \textsc{Klaassen, C. A. J.}, \textsc{Ritov,
Y.} and \textsc{Wellner, J. A.} (1993). \textit{Efficient and
Adaptive Estimation for Semiparametric Models}. Johns Hopkins
Univ. Press, Baltimore.
\MR{1245941}

\bibitem[\protect\citeauthoryear{Brookmeyer and Gail}{1994}]{r7}
\textsc{Brookmeyer, R.} and \textsc{Gail, M. H.} (1994). \textit{AIDS
Epidemiology: A Quantitative Approach}.
\textit{Monographs in Epidemiology and Biostatistics} \textbf{23}.
Oxford Univ. Press, New York.

\bibitem[\protect\citeauthoryear{Brunel}{2008}]{r8}
\textsc{Brunel, N.} (2008). Parameter estimation of ODE's via
nonparametric estimators. \textit{Electron. J. Statist.}
\textbf{2} 1242--1267.
\MR{2471285}

\bibitem[\protect\citeauthoryear{Burnham and Anderson}{2004}]{r9}
\textsc{Burnham, K. P.} and \textsc{Anderson, D. R.} (2004).
Multimodel inference: Understanding AIC and BIC in model selection.
\textit{Sociol. Methods Res.} \textbf{33} 261.
\MR{2086350}

\bibitem[\protect\citeauthoryear{Chappel and Godfrey}{1992}]{Chappel1992}
\textsc{Chappel, M. J.} and \textsc{Godfrey, K. R.} (1992). Structural
identifiability of the parameters of a nonlinear batch reactor model.
\textit{Math. Biosci.} \textbf{108} 245--251.
\MR{1154720}

\bibitem[\protect\citeauthoryear{Chen and Wu}{2008}]{r10}
\textsc{Chen, J.} and \textsc{Wu, H.} (2008). Efficient local estimation for
time-varying coefficients in
deterministic dynamic models with applications to HIV-1 dynamics.
\textit{J. Amer. Statist. Assoc.} \textbf{103} 369--384.
\MR{2420240}

\bibitem[\protect\citeauthoryear{Chen, He and Church}{1999}]{r11}
\textsc{Chen, T., He, H. L.} and \textsc{Church, G. M.} (1999). Modeling gene
expression with differential
equations. \textit{Pac. Symp. Biocomput.} 29--40.

\bibitem[\protect\citeauthoryear{Cobelli, Lepschy and Jacur}{1979}]{Cobelli1979}
\textsc{Cobelli, C.}, \textsc{Lepschy, A.} and \textsc{Jacur, R.} (1979).
Identifiability of compartmental systems and related structural
proerties. \textit{Math. Biosci.} \textbf{44} 1--18.

\bibitem[\protect\citeauthoryear{Daley and Gani}{1999}]{r12}
\textsc{Daley, D. J.} and \textsc{Gani, J.} (1999). \textit{Epidemic Modeling}.
Cambridge Univ. Press, Cambridge.

\bibitem[\protect\citeauthoryear{de Boor}{1978}]{r13}
\textsc{de Boor, C.} (1978). \textit{A Practical Guide to Splines}.
Springer, New York.
\MR{0507062}

\bibitem[\protect\citeauthoryear{Delgado}{1992}]{r14}
\textsc{Delgado, M. A.} (1992). Semiparametric generalized least
squares in the multivariate nonlinear regression model.
\textit{Econometric Theory} \textbf{8} 203--222.
\MR{1179510}

\bibitem[\protect\citeauthoryear{Donnet and Samson}{2007}]{r15}
\textsc{Donnet, S.} and \textsc{Samson, A.} (2007). Estimation of
parameters in incomplete data models defined by dynamical systems.
\textit{J. Statist. Plann. Inference} \textbf{137} 2815--2831.
\MR{2323793}

\bibitem[\protect\citeauthoryear{Englezos and Kalogerakis}{2001}]{r16}
\textsc{Englezos, P.} and \textsc{Kalogerakis, N.} (2001). \textit{Applied
Parameter Estimation for Chemical Engineers}.
Dekker, New York.

\bibitem[\protect\citeauthoryear{Grenander}{1981}]{r17}
\textsc{Grenander, U.} (1981). \textit{Abstract Inference}. Wiley,
New York.
\MR{0599175}

\bibitem[\protect\citeauthoryear{Hairer, N{\o}rsett and Wanner}{1993}]{r18}
\textsc{Hairer, E.}, \textsc{N{\o}rsett, S. P.} and \textsc{Wanner,
G.} (1993). \textit{Solving Ordinary Differential Equations I:
Nonstiff Problems}, 2nd ed. Springer, Berlin.
\MR{1227985}

\bibitem[\protect\citeauthoryear{He, Fung and Zhu}{2002}]{r19}
\textsc{He, X.}, \textsc{Fung, W. K.} and \textsc{Zhu, Z. Y.}
(2002). Estimation in a semiparametric model for longitudinal data
with unspecified dependence structure. \textit{Biometrika}
\textbf{89} 579--590.
\MR{1929164}

\bibitem[\protect\citeauthoryear{Heckman and Ramsay}{2000}]{r21}
\textsc{Heckman, N. E.} and \textsc{Ramsay, J. O.} (2000). Penalized
regression with model-based penalties.
\textit{Canad. J. Statist.} \textbf{28} 241--258.
\MR{1792049}

\bibitem[\protect\citeauthoryear{Ho et al.}{1995}]{r22}
\textsc{Ho, D. D.}, \textsc{Neumann, A. U.}, \textsc{Perelson, A.
S. et al.} (1995). Rapid turnover of plasma virions and CD4
lymphocytes in HIV-1 infection. \textit{Nature} \textbf{373}
123--126.

\bibitem[\protect\citeauthoryear{Huang}{1996}]{r23}
\textsc{Huang, J.} (1996). Efficient estimation for the proportinal
hazards model with interval censoring. \textit{Ann. Statist.}
\textbf{24} 540--568.
\MR{1394975}

\bibitem[\protect\citeauthoryear{Huang}{1999}]{r24}
\textsc{Huang, J.} (1999). Efficient estimation of the partly linear
additive Cox model. \textit{Ann. Statist.} \textbf{27} 1536--1563.
\MR{1742499}

\bibitem[\protect\citeauthoryear{Huang and Rossini}{1997}]{r25}
\textsc{Huang, J.} and \textsc{Rossini, A. J.} (1997). Sieve
estimation for the proportional-odds failure-time regression model
with interval censoring. \textit{J. Amer. Statist. Assoc.}
\textbf{92} 960--967.
\MR{1482126}

\bibitem[\protect\citeauthoryear{Huang}{2003}]{r26}
\textsc{Huang, J. Z.} (2003). Local asymptotics for polynomial
spline regression. \textit{Ann. Statist.} \textbf{31} 1600--1635.
\MR{2012827}

\bibitem[\protect\citeauthoryear{Huang, Zhang and Zhou}{2007}]{r27}
\textsc{Huang, J. Z.}, \textsc{Zhang, L.} and \textsc{Zhou, L.}
(2007). Efficient estimation in marginal partially linear models for
longitudinal/clustered data using plines.
\textit{Scand. J. Statist.} \textbf{34} 451--477.
\MR{2368793}

\bibitem[\protect\citeauthoryear{Huang, Liu and Wu}{2006}]{r28}
\textsc{Huang, Y.}, \textsc{Liu, D.} and \textsc{Wu, H.} (2006).
Hierarchical Bayesian methods for
estimation of parameters in a longitudinal HIV dynamic system.
\textit{Biometrics} \textbf{62}
413--423.
\MR{2227489}

\bibitem[\protect\citeauthoryear{Huang, Rosenkranz and Wu}{2003}]{HRW03}
\textsc{Huang, Y., Rosenkranz, S. L.} and \textsc{Wu, H.} (2003).
Modeling HIV dynamics and
antiviral response with consideration of time-varying drug exposures,
adherence and phenotypic sensitivity. \textit{Math. Biosci.}
\textbf{184} 165--186.

\bibitem[\protect\citeauthoryear{Huber}{1967}]{r29}
\textsc{Huber, P. J.} (1967). The behavior of maximum likelihood
estimates under nonstandard conditions. In \textit{Proc. Fifth Berkeley
Symp. Math. Statist. Probab.} 221--233.
Univ. California Press, Berkeley.
\MR{0216620}

\bibitem[\protect\citeauthoryear{Jansson and Revesz}{1975}]{r30}
\textsc{Jansson, B.} and \textsc{Revesz, L.} (1975).
Analysis of the growth of tumor cell populations. \textit{Math.
Biosci.} \textbf{19} 131--154.

\bibitem[\protect\citeauthoryear{Jeffrey and Xia}{2005}]{r31}
\textsc{Jeffrey, A. M.} and \textsc{Xia, X.} (2005). Identifiability of
HIV/AIDS model. In \textit{Deterministic and
Stochastic Models of AIDS Epidemics and HIV Infections with
Intervention} (W. Y. Tan and H.~Wu, eds.) 255--286. World Scientific, Singapore.

\bibitem[\protect\citeauthoryear{Jennerich}{1969}]{r32}
\textsc{Jennerich, R. I.} (1969). Asymptotic properties of
non-linear least squares estimators. \textit{Ann. Math. Statist.}
\textbf{40} 633--643.
\MR{0238419}

\bibitem[\protect\citeauthoryear{Jolliffe}{1972}]{Jolliffe1972}
\textsc{Jolliffe, I. T.} (1972). Discarding variables in a principal
component analyssis. I: Artificial data. \textit{J.~Roy. Statist. Soc. Ser. C Appl. Statist.}
\textbf{21} 160--172.
\MR{0311034}

\bibitem[\protect\citeauthoryear{Joshi, Seidel-Morgenstern and Kremling}{2006}]{r31a}
\textsc{Joshi, M., Seidel-Morgenstern, A.} and \textsc{Kremling, A.} (2006).
Exploiting the boostrap method
for quantifying parameter confidence intervals in dynamic systems.
\textit{Metabolic Engineering}
\textbf{8} 447--455.

\bibitem[\protect\citeauthoryear{Kolchin}{1973}]{Kolchin1973}
\textsc{Kolchin, E.} (1973). \textit{Differential Algebra and Algebraic
Groups}. Academic Press, New York.
\MR{0568864}

\bibitem[\protect\citeauthoryear{Kutta}{1901}]{r32a}
\textsc{Kutta, W.} (1901). Beitrag zur n\"{a}herungsweisen
itegration totaler differentialgleichungen. \textit{Zeitschr.
Math. Phys.} \textbf{46} 435--453.

\bibitem[\protect\citeauthoryear{Li et al.}{2002}]{r33}
\textsc{Li, L.}, \textsc{Brown, M. B.}, \textsc{Lee, K. H.} and
\textsc{Gupta, S.} (2002). Estimation and inference for a
spline-enhanced pupulation pharmacokinetic model.
\textit{Biometrics} \textbf{58} 601--611.
\MR{1933534}

\bibitem[\protect\citeauthoryear{Li, Osborne and Pravan}{2005}]{r34}
\textsc{Li, Z.}, \textsc{Osborne, M. R.} and \textsc{Pravan, T.}
(2005). Parameter estimation of ordinary differential equations.
\textit{IMA J. Numer. Anal.} \textbf{25} 264--285.
\MR{2126204}

\bibitem[\protect\citeauthoryear{Liang and Wu}{2008}]{r35}
\textsc{Liang, H.} and \textsc{Wu, H.} (2008). Parameter estimation for
differential equation models using a framework of measurement error
in regression models. \textit{J. Amer. Statist.
Assoc.} \textbf{103} 1570--1583.
\MR{2504205}

\bibitem[\protect\citeauthoryear{Ljung and Glad}{1994}]{Ljung1994}
\textsc{Ljung, L.} and \textsc{Glad, T.} (1994). On global
identifiability for arbitrary model parametrizations.
\textit{Automatica} \textbf{30} 265--276.
\MR{1261705}

\bibitem[\protect\citeauthoryear{Ma and Kosorok}{2005}]{r36}
\textsc{Ma, S.} and \textsc{Kosorok, M. R.} (2005). Robust
semiparametric M-estimation and the weighted bootstrap.
\textit{J. Multivariate Anal.} \textbf{96} 190--217.
\MR{2202406}

\bibitem[\protect\citeauthoryear{Malinvaud}{1970}]{r37}
\textsc{Malinvaud, E.} (1970). The consistancy of nonlinear
regressions. \textit{Ann. Math. Statist.} \textbf{41} 956--969.
\MR{0261754}

\bibitem[\protect\citeauthoryear{Mattheij and Molenaar}{2002}]{r38}
\textsc{Mattheij, R.} and \textsc{Molenaar, J.} (2002).
\textit{Ordinary Differential Equations in Theory and Practice}.
SIAM, Philadelphia.
\MR{1946758}

\bibitem[\protect\citeauthoryear{Miao et al.}{2008}]{r39}
\textsc{Miao, H., Dykes, C., Demeter, L. M., Cavenaugh, J., Park, S. Y.,
Perelson, A. S.} and \textsc{Wu, H.} (2008). Modeling and estimation of
kinetic parameters and replicative fitness of HIV-1 from
flow-cytometry-based growth competition experiments.
\textit{Bull. Math. Biol.} \textbf{70} 1749--1771.
\MR{2430325}

\bibitem[\protect\citeauthoryear{Michelson and Leith}{1997}]{r40}
\textsc{Michelson, S.} and \textsc{Leith, J. T.} (1997). Tumor heterogeneity and
growth control. In \textit{Tumor
Heterogeneity and Growth Control} (J. A. Adam and N. Bellomo, eds.) 295--326.
Birkh\"{a}user, Boston.

\bibitem[\protect\citeauthoryear{Moles, Banga and Keller}{2004}]{r41}
\textsc{Moles, C. G., Banga, J. R.} and \textsc{Keller, K.} (2004).
Solving nonconvex climate control problems: Pitfalls
and algorithm performances. \textit{Appl. Soft. Comput.} \textbf{5} 35--44.

\bibitem[\protect\citeauthoryear{Nowak and May}{2000}]{r42}
\textsc{Nowak, M. A.} and \textsc{May, R. M.} (2000). \textit{Virus
Dynamics: Mathematical Principles of Immunology and Virology}.
Oxford Univ. Press, Oxford.
\MR{2009143}

\bibitem[\protect\citeauthoryear{Ollivier}{1990}]{Ollivier1990}
\textsc{Ollivier, F.} (1990). Le probl\`{e}me de l'identifiabilit\'{e}
globale: \'{E}tude th\'{e} orique, m\'{e}thodes effectives et bornes de
complexit\'{e}. Ph.D. thesis, \'{E}cole Polytechnique, Paris, France.

\bibitem[\protect\citeauthoryear{Pakes and Pollard}{1989}]{r43}
\textsc{Pakes, A.} and \textsc{Pollard, D.} (1989). Simulation and
the asymptotics of optimization estimators. \textit{Econometrica}
\textbf{57} 1027--1057.
\MR{1014540}

\bibitem[\protect\citeauthoryear{Perelson et al.}{1997}]{r44}
\textsc{Perelson, A. S.}, \textsc{Essunger, P.}, \textsc{Cao, Y.},
\textsc{Vesanen, M.}, \textsc{Hurley, A.}, \textsc{Saksela, K.},
\textsc{Markowitz, M.} and  \textsc{Ho, D. D.} (1997). Decay characteristics of
HIV-1-infected compartments during combination therapy.
\textit{Nature} \textbf{387} 188--191.

\bibitem[\protect\citeauthoryear{Perelson and Nelson}{1999}]{r45}
\textsc{Perelson, A. S.} and \textsc{Nelson, P. W.} (1999). Mathematical
analysis of HIV-1 dynamics in vivo. \textit{SIAM Rev.}
\textbf{41} 3--44.
\MR{1669741}

\bibitem[\protect\citeauthoryear{Perelson et al.}{1996}]{r46}
\textsc{Perelson, A. S.}, \textsc{Neumann, A. U.},
\textsc{Markowitz, M.}, \textsc{Leonard, J. M.} and
\textsc{Ho, D. D.} (1996). HIV-1 dynamics
in vivo: Virion clearance rate, infected cell life-span, and viral
generation time. \textit{Science} \textbf{271} 1582--1586.

\bibitem[\protect\citeauthoryear{Pohjanpalo}{1978}]{Pohjanpalo1978}
\textsc{Pohjanpalo, H.} (1978). System identifiability based on the
power series expansion of the solution. \textit{Math. Biosci.}
\textbf{41} 21--33.
\MR{0507373}

\bibitem[\protect\citeauthoryear{Pollard}{1984}]{r47}
\textsc{Pollard, D.} (1984). \textit{Convergence of Stochastic
Processes}. Springer, New York.
\MR{0762984}

\bibitem[\protect\citeauthoryear{Pollard}{1985}]{r48}
\textsc{Pollard, D.} (1985). New ways to prove central limit
theorems. \textit{Econometric Theory} \textbf{1} 295--314.

\bibitem[\protect\citeauthoryear{Pollard}{1990}]{r49}
\textsc{Pollard, D.} (1990). \textit{Empirical Processes Theory and
Applications}. IMS, Hayward, CA.
\MR{1089429}

\bibitem[\protect\citeauthoryear{Poyton et al.}{2006}]{r50}
\textsc{Poyton, A. A.}, \textsc{Varziri, M. S.}, \textsc{McAuley, K.
B.}, \textsc{McLellen, P. J.} and \textsc{Ramsay, J. O.} (2006).
Parameter estimation in continuous-time dynamic models using
principal differential analysis. \textit{Computers and Chemical
Engineering} \textbf{30} 698--708.

\bibitem[\protect\citeauthoryear{Putter et al.}{2002}]{r51}
\textsc{Putter, H.}, \textsc{Heisterkamp, S. H.}, \textsc{Lange, J.
M.} and \textsc{de Wolf, F.} (2002). A Bayesian approach to
parameter estimation in HIV dynamical models. \textit{Stat. Med.}
\textbf{21} 2199--2214.

\bibitem[\protect\citeauthoryear{Quaiser and M\"{o}nnigmann}{2009}]{Quaiser2009}
\textsc{Quaiser, T.} and \textsc{M\"{o}nnigmann, M.} (2009). Systematic
identifiability testing for unambiguous mechanistic modeling---application
to JAK-STAT, MAP kinase, and NF-$\kappa$B signaling pathway
models. \textit{BMC Sys. Bio.} \textbf{3} 50.

\bibitem[\protect\citeauthoryear{Ramsay}{1996}]{r52}
\textsc{Ramsay, J. O.} (1996). Principal Differential Analysis: Data
Reduction by Differential Operators. \textit{J.~Roy. Statist. Soc. Ser. B}
\textbf{58} 495--508.
\MR{1394362}

\bibitem[\protect\citeauthoryear{Ramsay et al.}{2007}]{r53}
\textsc{Ramsay, J. O.}, \textsc{Hooker, G.}, \textsc{Campbell, D.} and
\textsc{Cao, J.} (2007).
Parameter estimation for differential equations: A generalized
smoothing approach (with discussion).
\textit{J. R. Stat. Soc. Ser. B Stat. Methodol.} \textbf{69} 741--796.
\MR{2368570}

\bibitem[\protect\citeauthoryear{Ritt}{1950}]{Ritt1950}
\textsc{Ritt, J. F.} (1950). \textit{Differential Algebra}.
Amer. Math. Soc., Providence, RI.
\MR{0035763}

\bibitem[\protect\citeauthoryear{Rodriguez-Fernandez, Egea and
Banga}{2006}]{r54}
\textsc{Rodriguez-Fernandez, M.}, \textsc{Egea, J. A.} and
\textsc{Banga, J. R.} (2006). Novel metaheuristic for parameter
estimation in nonlinear dynamic biological systems.
\textit{BMC Bioinformatics} \textbf{7} 1--18.

\bibitem[\protect\citeauthoryear{Runge}{1895}]{r55}
\textsc{Runge, C.} (1895). Ueber die numerische Aufl\"{o}sung von
Differentialgleichungen. \textit{Math. Ann.} \textbf{46} 167--178.
\MR{1510879}

\bibitem[\protect\citeauthoryear{Schick}{1986}]{r56}
\textsc{Schick, A.} (1986). On asymptotically efficient estimation
in semiparametric models. \textit{Ann. Statist.} \textbf{14}
1139--1151.
\MR{0856811}

\bibitem[\protect\citeauthoryear{Schumaker}{1981}]{r57}
\textsc{Schumaker, L. L.} (1981). \textit{Spline Functions}. Wiley,
New York.
\MR{0606200}

\bibitem[\protect\citeauthoryear{Shen}{1997}]{r58}
\textsc{Shen, X.} (1997). On methods of sieves and penalization.
\textit{Ann. Statist.} \textbf{25} 2555--2591.
\MR{1604416}

\bibitem[\protect\citeauthoryear{Shen and Wong}{1994}]{r59}
\textsc{Shen, X.} and \textsc{Wong, W. H.} (1994). Convergence rate
of sieve estimates. \textit{Ann. Statist.} \textbf{22} 580--615.
\MR{1292531}

\bibitem[\protect\citeauthoryear{Stone}{1982}]{r60}
\textsc{Stone, C. J.} (1982). Optimal global rates of convergence
for nonparametric regression. \textit{Ann. Statist.} \textbf{10}
1040--1053.
\MR{0673642}

\bibitem[\protect\citeauthoryear{Stone}{1985}]{r61}
\textsc{Stone, C. J.} (1985). Additive regression and other
nonparametric models. \textit{Ann. Statist.} \textbf{13} 689--705.
\MR{0790566}

\bibitem[\protect\citeauthoryear{Storn and Price}{1997}]{r62}
\textsc{Storn, R.} and \textsc{Price, K.} (1997). Differential
evolution---a simple and efficient heuristic for
global optimization over continuous spaces.
\textit{J. Global Optim.} \textbf{11} 341--359.
\MR{1479553}

\bibitem[\protect\citeauthoryear{Swartz and Bremermann}{1975}]{r63}
\textsc{Swartz, J.} and \textsc{Bremermann, H.} (1975). Discussion of
parameter estimation in
biological modeling: Algorithms for estimation and evaluation of the
estimates. \textit{J. Math. Biol.} \textbf{1} 241--275.

\bibitem[\protect\citeauthoryear{Tan and Wu}{2005}]{r64}
\textsc{Tan, W. Y.} and \textsc{Wu, H.} (2005). \textit{Deterministic and
Stochastic Models of AIDS
Epidemics and HIV Infections With Intervention}. World
Scientific, Singapore.
\MR{2169300}

\bibitem[\protect\citeauthoryear{Thomaseth et al.}{1996}]{r65}
\textsc{Thomaseth, K.}, \textsc{Alexandra, K. W.}, \textsc{Bernhard,
L. et al.} (1996). Integrated mathematical model to assess
$\beta$-cell activity during the oral glucose test. \textit{Amer. J.
Phisiol.} \textbf{270} E522--E531.

\bibitem[\protect\citeauthoryear{Vajda et al.}{1989}]{Vajada1989}
\textsc{Vajda, S., Rabitz, H., Walter, E.} and \textsc{Lecourtier, Y.} (1989).
Qualitative and quantitative identifiability analysis of nonlinear
chemical kinetiv-models. \textit{Chem. Eng. Commun.} \textbf{83} 191--219.

\bibitem[\protect\citeauthoryear{van der Vaart and Wellner}{1996}]{r66}
\textsc{van der Vaart, A. W.} and \textsc{Wellner, J. A.} (1996).
\textit{Weak Convergence and Empirical Processes}. Springer,
New York.
\MR{1385671}

\bibitem[\protect\citeauthoryear{van Domselaar and Hemker}{1975}]{r67}
\textsc{van Domselaar, B.} and \textsc{Hemker, P. W.} (1975).
Nonlinear parameter estimation in initial value problmes.
Report NW18/75, Math. Centrum, Amsterdam.

\bibitem[\protect\citeauthoryear{Varah}{1982}]{r68}
\textsc{Varah, J. M.} (1982). A spline least squares method for
numerical parameter estimation in differential equations.
\textit{SIAM J. Sci. Comput.} \textbf{3} 28--46.
\MR{0651865}

\bibitem[\protect\citeauthoryear{Varziri et al.}{2008}]{r69}
\textsc{Varziri, M. S.}, \textsc{Poyton, A. A.}, \textsc{McAuley, K.
B.}, \textsc{McLellen, P. J.}
and \textsc{Ramsay, J. O.} (2008). Selecting optimal weighting factors
in iPDA for parameter
estimation in continuous-time dynamic models.
\textit{Comp. Chem. Eng.} \textbf{32} 3011--3022.

\bibitem[\protect\citeauthoryear{Walter}{1987}]{Walter1987}
\textsc{Walter, E.} (1987). \textit{Identifiability of Parameteric
Models}. Pergamon Press, Oxford.

\bibitem[\protect\citeauthoryear{Wei et al.}{1995}]{r70}
\textsc{Wei, X.}, \textsc{Ghosh, S. K.}, \textsc{Taylor, M. E. et
al.} (1995). Viral dynamics in human immunodeficiency virus type 1
infection. \textit{Nature} \textbf{373} 117--122.

\bibitem[\protect\citeauthoryear{Wellner and Zhang}{2007}]{r71}
\textsc{Wellner, J. A.} and \textsc{Zhang, Y.} (2007). Two
likelihood-based semiparametric estimation methods for panel count
data with covariates. \textit{Ann. Statist.} \textbf{35}
2106--2142.
\MR{2363965}

\bibitem[\protect\citeauthoryear{Wu}{1981}]{r72}
\textsc{Wu, C. F.} (1981). Asymptotic theory of nonlinear least
squares estiamtion. \textit{Ann. Statist.} \textbf{9} 501--513.
\MR{0615427}

\bibitem[\protect\citeauthoryear{Wu}{2005}]{r73}
\textsc{Wu, H.} (2005). Statistical methods for HIV dynamic studies
in AIDS clinical trials. \textit{Stat. Methods Med. Res.}
\textbf{14} 1--22.
\MR{2135921}

\bibitem[\protect\citeauthoryear{Wu et al.}{2005}]{r74}
\textsc{Wu, H.}, \textsc{Huang, Y.}, \textsc{Acosta, E. P. et al.}
(2005). Modeling long-term HIV dynamics and antiretroviral response:
Effects of drug potency, pharmacokinetics, adherence, and drug
resistance. \textit{JAIDS} \textbf{39} 272--283.

\bibitem[\protect\citeauthoryear{Wu et al.}{2008}]{Wu2008}
\textsc{Wu, H., Zhu, H., Miao, H.} and \textsc{Perelson, A. S.} (2008).
Identifiability and
statistical estimation of dynamic parameters in HIV/AIDS dynamic
models. \textit{Bull. Math. Biol.}
\textbf{70} 785--799.
\MR{2393024}

\bibitem[\protect\citeauthoryear{Xia}{2003}]{Xia2003}
\textsc{Xia, X.} (2003). Estimation of HIV/AIDS parameters.
\textit{Automatica J. IFAC} \textbf{39} 1983--1988.
\MR{2142834}

\bibitem[\protect\citeauthoryear{Xia and Moog}{2003}]{XiaMoog2003}
\textsc{Xia, X.} and \textsc{Moog, C. H.} (2003). Identifiability of
nonlinear systems with
applications to HIV/AIDS models. \textit{IEEE Trans. Automat. Control}
\textbf{48} 330--336.
\MR{1957979}

\bibitem[\protect\citeauthoryear{Xue, Lam and Li}{2004}]{r76}
\textsc{Xue, H.}, \textsc{Lam, K. F.} and \textsc{Li, G.} (2004).
Sieve maximum likelihood estimator for semiparametric regression
models with current status data. \textit{J. Amer. Statist. Assoc.}
\textbf{99} 346--356.
\MR{2062821}

\end{thebibliography}
\end{document}